\documentclass[10pt,oneside,a4paper,reqno]{amsart}

\usepackage[T1]{fontenc}
\usepackage{lmodern}
\usepackage{amsmath,amssymb,amsthm,bm,comment}

\usepackage[
    a4paper,
    margin=1.1in,
    top=1.2in,
    bottom=1in,
    headheight=15pt,
    headsep=0.4in,
    footskip=1.2cm
]{geometry}

\usepackage[
    activate={true,nocompatibility},
    final,
    tracking=true,
    spacing=true,
    stretch=15,
    shrink=15
]{microtype}

\microtypecontext{spacing=nonfrench}
\SetTracking{encoding=*,family=*}{-20}

\usepackage{etoolbox}
\usepackage[dvipsnames]{xcolor}

\definecolor{darkorange}{HTML}{85423e}
\definecolor{domainblue}{HTML}{294E66}
\definecolor{fillteal}{HTML}{387D76}
\definecolor{coregold}{HTML}{A06A24}

\usepackage{graphicx}
\usepackage{tikz}
\usetikzlibrary{arrows.meta,calc,patterns}

\usepackage{caption}
\usepackage[
    colorlinks=true,
    linkcolor=darkorange,
    citecolor=darkorange,
    urlcolor=darkorange
]{hyperref}

\numberwithin{equation}{section}

\theoremstyle{plain}
\newtheorem{theorem}{Theorem}
\newtheorem{proposition}[theorem]{Proposition}
\newtheorem{lemma}[theorem]{Lemma}

\theoremstyle{remark}
\newtheorem{remark}{Remark}

\AtBeginEnvironment{remark}{\pushQED{\qed}}
\AtEndEnvironment{remark}{\popQED}

\newcommand{\R}{\mathbb{R}}
\newcommand{\C}{\mathbb{C}}
\newcommand{\e}{\varepsilon}
\newcommand{\Sone}{\mathbb{S}^1}
\newcommand{\TT}{\mathbb{T}}
\newcommand{\cA}{\mathcal{A}}
\newcommand{\cX}{\mathcal{X}}

\newcommand{\bn}{\bm{n}}
\newcommand{\bu}{\bm{u}}
\newcommand{\bU}{\bm{U}}
\newcommand{\bv}{\bm{v}}
\newcommand{\bV}{\bm{V}}
\newcommand{\bw}{\bm{w}}
\newcommand{\bW}{\bm{W}}
\newcommand{\bh}{\bm{h}}
\newcommand{\bH}{\bm{H}}
\newcommand{\bQ}{\bm{Q}}
\newcommand{\bP}{\bm{P}}
\newcommand{\bz}{\bm{z}}
\newcommand{\bxi}{\bm{\xi}}
\newcommand{\bZ}{\bm{Z}}
\newcommand{\btau}{\bm{\tau}}
\newcommand{\Ucomp}{U}

\newcommand{\E}{\mathcal{E}}
\newcommand{\cF}{\mathcal{F}}
\newcommand{\cS}{\mathcal{S}}
\newcommand{\Surf}{S}
\newcommand{\J}{J}
\newcommand{\mdeg}{m}
\newcommand{\Bp}{\mathcal{B}_p}
\newcommand{\cp}{c_p}
\newcommand{\lc}{\ell}
\newcommand{\rh}{\rho}
\newcommand{\ps}{\psi}
\newcommand{\vp}{\varphi}
\newcommand{\am}{\alpha}
\newcommand{\ad}{\alpha}
\newcommand{\qcut}{q}
\newcommand{\qexp}{q}
\newcommand{\Dp}{D}
\newcommand{\Zset}{Z}
\newcommand{\Dext}{D}
\newcommand{\Ddef}{D}
\newcommand{\Herr}{H}
\newcommand{\Hmor}{H}
\newcommand{\rpack}{\rho}
\newcommand{\dd}{\,dx}

\title{Asymptotic radial symmetry of $p$-Ginzburg--Landau minimizers on the disk}
\date{}

\author{Giovanni Di Fratta}
\address{Giovanni Di Fratta: Dipartimento di Matematica e Applicazioni ``R. Caccioppoli'', Universit\`a degli Studi di Napoli ``Federico II'', Via Cintia, 80126 Napoli, Italy.}
\email{giovanni.difratta@unina.it}

\author{Valeriy V. Slastikov}
\address{Valeriy V. Slastikov: School of Mathematics, University of Bristol, Bristol BS8 1TW, United Kingdom.}
\email{valeriy.slastikov@bristol.ac.uk}

\author{Arghir D. Zarnescu}
\address{Arghir D. Zarnescu: BCAM, Basque Center for Applied Mathematics, Mazarredo 14, E48009 Bilbao, Bizkaia, Spain; IKERBASQUE, Basque Foundation for Science, Maria Diaz de Haro 3, 48013 Bilbao, Bizkaia, Spain; Simion Stoilow Institute of Mathematics of the Romanian Academy, P.O. Box 1-764, 014700 Bucharest, Romania.}
\email{azarnescu@bcamath.org}



\begin{document}

\begin{abstract}
We solve the minimizer case of Brezis's Open
Problem~2.3 for every fixed $p>2$.
For $\e>0$ and 
$\bu_\e\in W^{1,p}(B_1;\R^2)$ a global
minimizer of the $p$-Ginzburg--Landau energy
\[
 \int_{B_1}\left(
 \frac{|\nabla\bu|^p}{p}
 +\frac{(1-|\bu|^2)^2}{4\e^2}
 \right)\,dx
\]
with prescribed trace
$(\bu_\e)_{|\partial B_1}(x)=x$,
we prove that, as $\e\downarrow0$,
$\bu_\e(x)\to x/|x|$ at every $x$ in the punctured unit disk 
$B_1\setminus\{0\}$.
The convergence is strong in
$W^{1,p}_{\mathrm{loc}}(B_1\setminus\{0\};\R^2)$
and hence uniform on every compact subset
of the punctured disk.
On each such compact subset, the potential energy
tends to zero and the total energy converges to the
$p$-Dirichlet energy of the radial map.
No symmetry or uniqueness of the minimizers
is assumed.
The proof combines whole-plane concentration,
a strict exterior filling estimate, and quantitative
convexity to locate the concentration point and
identify the limiting map.
\end{abstract}

\maketitle

\section{The problem and the distinction from the quadratic case}
For $a\in\R^2$ and $r\in(0,\infty)$, let
$B_r(a)=\{x\in\R^2:|x-a|<r\}$ and $B_r=B_r(0)$.
We identify $\R^2$ with $\C$ when multiplying maps and write
$\Sone=\{z\in\C:|z|=1\}$ and $\TT=\R/(2\pi\mathbb Z)$. Fix $p\in(2,\infty)$ and  set
\begin{equation}\label{eq:class}
 \cA=\{\bu\in W^{1,p}(B_1;\R^2):\bu_{|\partial B_1}(x)=x
                  \text{ on }\partial B_1\}.
\end{equation}
For $\e>0$, an open set $G\subset \R^2$, a map
$\bu\in W^{1,p}_{\mathrm{loc}}(G;\R^2)$, and a measurable set $D\subset G$, define
\begin{equation}\label{eq:energy}
 \E_\e(\bu;D)=\int_D\left(\frac{|\nabla \bu|^p}{p}
               +\frac{(1-|\bu|^2)^2}{4\e^2}\right)\dd,
 \qquad \E_\e(\bu)=\E_\e(\bu;B_1)\quad\text{for }\bu\in\cA.
\end{equation}
We use the Frobenius norm for matrices in $\R^{2\times2}$,
in particular when writing $|\nabla\bu|$.
For every $\e>0$, the nonnegative functional
$\E_\e$ is well defined and finite on $\cA$.
For general $G$, $D$, and $\bu$ as above,
$\E_\e(\bu;D)$ may equal $+\infty$; for example,
$\E_\e(0;\R^2)=+\infty$.
Since $p>2$, maps in $W^{1,p}_{\mathrm{loc}}(G;\R^2)$ have
continuous representatives on $G$, which we always use.
For maps in $W^{1,p}$ on a bounded Lipschitz domain, these
representatives extend continuously to its closure. Throughout, in estimates, $C,C_p$ denote strictly positive finite constants
depending only on $p$ and may increase between estimates;
additional dependencies are indicated by subscripts.
The boundary condition singles out the radial map
\begin{equation}\label{eq:radial}
 \bn\in C^\infty(\R^2\setminus\{0\};\Sone),
 \qquad \bn(x)=\frac{x}{|x|}.
\end{equation}
Although $\bn\notin W^{1,p}(B_1;\R^2)$, its $p$-Dirichlet energy
outside any disk centered at zero is finite. Indeed, $|\nabla \bn(x)|=|x|^{-1}$,
and polar integration gives
\begin{equation}\label{eq:rad_tail}
 \frac1p\int_{\R^2\setminus B_r}|\nabla \bn|^p\,dx
 =\frac{2\pi}{p}\int_r^\infty t^{1-p}\,dt
 =\cp r^{2-p},
 \qquad
 \cp=\frac{2\pi}{p(p-2)}>0,
 \qquad r>0.
\end{equation}
In particular,
$\cp\overset{\eqref{eq:rad_tail}}{=}
p^{-1}\int_{\R^2\setminus B_1}|\nabla \bn|^p\,dx$:
$\cp$ is precisely the $p$-Dirichlet energy of the radial map $\bn$ outside the unit disk.

Brezis asks in \cite[Open Problem~2.3]{Brezis}, also formulated
in \cite[Open Problem~30]{BM}, whether every family
$(\bu_\e)_{\e>0}\subset\cA$ of global minimizers of $\E_\e$ satisfies
$\bu_\e(x)\to \bn(x)$ for every $x\in B_1\setminus\{0\}$ as
$\e\downarrow0$, and whether the same holds for arbitrary
critical points. We give an affirmative answer to the first
question, concerning global minimizers, by proving something
even stronger: uniform convergence on $\overline{B_1\setminus B_r}$
and strong convergence in $W^{1,p}(B_1\setminus\overline{B_r};\R^2)$
for every fixed $r\in(0,1)$.
\begin{theorem}\label{thm:main}
For each $\e>0$, let $\bu_\e\in\cA$ be a global minimizer of $\E_\e$.
For each fixed $r\in(0,1)$, set
$A_r=B_1\setminus\overline{B_r}$. Then
\begin{equation}\label{eq:conv}
 \bu_\e\longrightarrow \bn
 \quad\text{uniformly on }\overline{A_r}
 \quad\text{and strongly in }W^{1,p}(A_r;\R^2)
 \qquad\text{as }\e\downarrow0.
\end{equation}
Moreover,
\begin{equation}\label{eq:en_lim}
 \frac1{4\e^2}\int_{A_r}(1-|\bu_\e|^2)^2\dd\longrightarrow0,
 \qquad
 \E_\e(\bu_\e;A_r)\longrightarrow \cp(r^{2-p}-1)
 \qquad\text{as }\e\downarrow0.
\end{equation}
Thus, on every fixed annulus away from the origin, the potential
energy vanishes, and the total energy converges to the
$p$-Dirichlet energy of the radial map $\bn$.
\end{theorem}

The proof of Theorem~\ref{thm:main} will be given in Section~\ref{sec:proof}
once we have established the necessary building blocks in the next sections.

\begin{remark}[Local convergence and convergence up to the boundary]
Theorem~\ref{thm:main} implies
\begin{equation}
 \bu_\e\longrightarrow\bn
 \quad\text{strongly in }
 W^{1,p}_{\mathrm{loc}}(B_1\setminus\{0\};\R^2)
 \qquad\text{as }\e\downarrow0.
\end{equation}
Since $p>2$, Sobolev embedding also gives uniform
convergence on every compact subset of the punctured
disk.
The theorem provides additional control near the
outer boundary. For every fixed $r\in(0,1)$,
the convergence is strong in
$W^{1,p}(B_1\setminus\overline{B_r};\R^2)$
and therefore uniform on the closed annulus
$\{x\in\R^2:r\leqslant|x|\leqslant1\}$.
Thus convergence holds up to $\partial B_1$:
uniformly on every compact subset of
$\overline{B_1}\setminus\{0\}$, with strong
$W^{1,p}$ control on each fixed annulus reaching
the outer boundary.
\end{remark}

\begin{remark}[Why the proof differs from the quadratic case]
For $p=2$, the asymptotic convergence considered here
is already known for minimizers \cite{BBH} and for
critical points \cite{PR}, as recalled immediately
after \cite[Open Problem~2.3]{Brezis}.
The passage to $p>2$ introduces the quasilinear
difficulties described above and removes the exact
quadratic structure of the Dirichlet energy.
In particular, the modulus and phase contributions
are coupled through the $p$-power of the gradient,
so their separate control requires nonlinear
estimates that retain enough information to identify
the limiting map.
Our proof addresses them through two quantitative
estimates: the strict exterior improvement in
Lemma~\ref{lem:fill}, which yields decay of the
exterior energy excess, and the convexity remainder
in Proposition~\ref{prop:defect}, which converts
this decay into control of the modulus and phase.
Thus the argument requires more than compactness.

The regime $p>2$ nevertheless supplies two properties
that make this strategy possible: Morrey compactness
on bounded Lipschitz domains in $\R^2$ and the finite
radial exterior energy
$\cp r^{2-p}\to0$ as $r\to\infty$.
At $p=2$, the embedding of $W^{1,2}$ into $C^0$
fails, while
\begin{equation}\label{eq:ann_en}
 \frac12\int_{B_R\setminus B_r}
 |\nabla\bn|^2\dd
 =\pi\log(R/r),
 \qquad 0<r<R<\infty.
\end{equation}
The radial exterior energy is therefore infinite.
These differences concern the foundations of our
comparison argument: the quadratic case cannot
be recovered by letting $p\downarrow2$ in this proof.

The finite whole-plane energy available for $p>2$
also allows us to use the results of
Almog--Berlyand--Golovaty--Shafrir on degree reduction
and concentration \cite{ABGS}.
Their degree-reduction and
compactness results allow us to confine the modulus defects
of arbitrary disk minimizers to a single small
region. Determining its position and the limiting
map requires additional information from the
prescribed boundary values, which enters through
the exterior comparison and convexity estimates
described above. The radial analysis in
\cite{ABGSrad} does not resolve these questions
without an a priori symmetry assumption.
Lei obtained related results on singular sets,
potential decay, and energy concentration
\cite{Lei07,Lei10,Lei22}, but our proof is
independent of those estimates.
\end{remark}

The remainder of the paper is organized as follows. In 
Section~\ref{sec:conc}, we show that each fixed
modulus-defect  set is confined to a ball of
radius $O(\lc)$, where $\lc=\e^{2/p}$.
To control the energy outside this ball,
Section~\ref{sec:fill} constructs an annular
interpolation from its possibly nonradial boundary
trace to the prescribed radial map.
The strict gain in this construction yields decay
of the exterior energy excess in
Section~\ref{sec:excess}.
Section~\ref{sec:defect} then converts
this energy bound into estimates for the modulus
gradient, phase gradient, and potential energy
of the extended minimizer.
Finally, Section~\ref{sec:proof} uses these
estimates and the prescribed exterior map to prove
$a_\e\to0$ and the convergence of $\bu_\e$
to $\bn$ on every fixed annulus away
from the origin.

\section{The natural scale and concentration into one core}
\label{sec:conc}

For every $\e>0$, the direct method in the calculus of variations
yields a global minimizer $\bu_\e\in\cA$ of $\E_\e$.
Every global minimizer satisfies $|\bu_\e|\leqslant1$ on
$\overline{B_1}$. Indeed, the nearest-point projection onto the
closed unit disk defined by
\begin{equation}\label{eq:proj}
 \bP:\R^2\to\overline{B_1},
 \qquad \bP(z)=\frac{z}{\max\{1,|z|\}},
\end{equation}
fixes $\overline{B_1}$, retracts $\R^2\setminus B_1$ to $\mathbb{S}^1$, and is $1$-Lipschitz\footnote{For
$z\in\R^2$ and $|q|\leqslant1$, the formula for $\bP$
gives $(z-\bP(z))\cdot(q-\bP(z))\leqslant0$.
Indeed, this is immediate if also $|z|\leqslant1$ because the LHS vanishes; otherwise,
$z-\bP(z)=(|z|-1)\bP(z)$ and $\bP(z)\cdot q\leqslant1$ by Cauchy–Schwarz given that 
$|\bP(z)|\leqslant1$ and $|q|\leqslant1$.
Applying this inequality with $(z,q)=(z,\bP(w))$ and
$(w,\bP(z))$ and adding yields
$|\bP(z)-\bP(w)|^2\leqslant(z-w)\cdot(\bP(z)-\bP(w))
\leqslant|z-w|\,|\bP(z)-\bP(w)|$.
Thus $|\bP(z)-\bP(w)|\leqslant|z-w|$ for all $z,w\in\R^2$.}. 
%
%
Hence the competitor
$\bP\circ \bu_\e\in\cA$ has no larger
Dirichlet energy. If $|\bu_\e(x_0)|>1$ at some $x_0\in B_1$,
continuity gives an open neighborhood of $x_0$ on which
$|\bu_\e|>1$. At every point where $|\bu_\e|>1$, the
projected map has modulus one, so its potential density is zero;
where $|\bu_\e|\leqslant1$,
the map is unchanged.
The potential energy therefore decreases strictly.
Consequently,
\begin{equation}\label{eq:proj_dec}
 \E_\e(\bP\circ \bu_\e)<\E_\e(\bu_\e),
\end{equation}
contradicting the minimality of $\bu_\e$. Thus
$|\bu_\e|\leqslant1$ throughout $\overline{B_1}$.

Every global minimizer $\bu_\e$ has at least one zero
$a_\e\in B_1$. Indeed, otherwise
$\bu_\e/|\bu_\e|$ would define a continuous map
$\overline{B_1}\to\Sone$ whose restriction to
$\partial B_1$ has degree one, contradicting the fact
that the boundary restriction of any such map has
degree zero. 
In particular, $\inf_{B_1}|\bu_\e|=0$.
\begin{remark}[Localization and isolation of zeros]
To the best of our knowledge, it is not known whether
the zeros of arbitrary global minimizers in the present
setting are isolated. For the radial whole-plane
minimizer, Almog--Berlyand--Golovaty--Shafrir prove
positivity of the radial profile away from the origin,
which is therefore its unique zero \cite{ABGSrad}.
Likewise, proving that every disk minimizer is
equivariant, together with positivity of its minimizing
radial profile, would establish uniqueness of its zero
at the origin.
Our argument requires neither isolation nor uniqueness.
The zero set
\begin{equation}
 \Zset_\e=\{x\in B_1:\bu_\e(x)=0\}
\end{equation}
is nonempty and compact, by continuity of $\bu_\e$ and
its nonvanishing boundary values.
Proposition~\ref{prop:conc} below, applied with
$\eta=1/2$, will show that, for sufficiently small $\e$,
\begin{equation}
 \Zset_\e\subset B_{C_{p,1/2}\lc}(a_\e),
 \qquad \lc=\e^{2/p},
\end{equation}
where $a_\e$ is any chosen zero and $C_{p,1/2}>0$
is independent of $\e$. This localizes the entire
zero set on a shrinking length scale, but does not
determine its cardinality or exclude accumulation
of zeros within that ball.
\end{remark}

To investigate the behavior of $\bu_\e$ near its zeros,
we introduce a rescaling that eliminates $\e$ from
the energy density and reduces the analysis to a
fixed functional on expanding domains. For that, we observe that
for any $\lc>0$, the rescaled map
$\bw(y)=\bu_\e(\lc y)$ belongs to
$W^{1,p}(B_{\lc^{-1}};\R^2)$, and $x=\lc y$ gives
\begin{equation}\label{eq:scale_gen}
 \E_\e(\bu_\e)
 =\lc^{2-p}\int_{B_{\lc^{-1}}}
 \left(\frac{|\nabla \bw|^p}{p}
 +\frac{\lc^p}{4\e^2}(1-|\bw|^2)^2\right)\,dy.
\end{equation}
Thus the choice $\lc=\e^{2/p}$ removes the parameter
from the integrand.
With this choice, set $R=\lc^{-1}>0$. The rescaling $\bw(y)=\bu(\lc y)$
is a bijection between $\cA$ and
\begin{equation}\label{eq:class_r}
 \cA_R=\{\bv\in W^{1,p}(B_R;\R^2):\bv_{|\partial B_R}=\bn\},
\end{equation}
and satisfies (here, $D\subset\R^2$ is open and
$\bv\in W^{1,p}_{\mathrm{loc}}(D;\R^2)$)
\begin{equation}\label{eq:scale}
 \E_\e(\bu)=\lc^{2-p}\cF(\bw;B_R),
 \qquad
 \cF(\bv;D)=\int_D\left(
 \frac{|\nabla \bv|^p}{p}+\frac{(1-|\bv|^2)^2}{4}
 \right)\dd.
\end{equation}
We define $\cF(\bv;D)$ for every open set $D\subset\R^2$
and every $\bv\in W^{1,p}_{\mathrm{loc}}(D;\R^2)$,
allowing the value $+\infty$. For $\bv\in\cA_R$, however,
the energy is finite: the gradient is $p$-integrable,
and Sobolev embedding makes $\bv$ bounded on $B_R$,
so the potential is integrable as well.

It is clear that a map $\bu\in\cA$ minimizes $\E_\e$
if and only if its rescaling $\bw\in\cA_R$, with $R=\lc^{-1}=\e^{-2/p}$, minimizes
$\cF(\cdot;B_R)$. Through the identity $\bu(x)=\bw(Rx)$,
the asymptotic behavior of minimizers as $\e\downarrow0$
is therefore equivalent to that of minimizers of $\cF$
on expanding disks as $R=\e^{-2/p}\to\infty$.
Under this rescaling, each fixed annulus
$A_r=B_1\setminus\overline{B}_r$, $r\in(0,1)$,
corresponds to $B_R\setminus\overline{B}_{rR}$.

To study this large-domain limit, we compare the minimum
energies on disks of different radii. For $R>0$, let
\begin{equation}\label{eq:j}
 \J(R):=\min_{\bv\in\cA_R}\cF(\bv;B_R)+\cp R^{2-p}.
\end{equation}
The term $\cp R^{2-p}$ in \eqref{eq:j} is exactly
the energy of $\bn$ outside $B_R$.
Thus $\J(R)$ is the whole-plane energy obtained by extending
a disk minimizer by its prescribed radial boundary map.
This interpretation also gives a monotonicity property.
For $R'>R>0$ and a minimizer $\bw_R\in\cA_R$, define
\begin{equation}\label{eq:ext_r}
 \widetilde \bw_{R,R'}(y)=
 \begin{cases}
  \bw_R(y),&y\in B_R,\\
  \bn(y),&y\in B_{R'}\setminus B_R.
 \end{cases}
\end{equation}
The traces agree on $\partial B_R$, so
$\widetilde \bw_{R,R'}\in\cA_{R'}$. Since $|\bn|=1$,
the extension adds only the radial Dirichlet energy:
\begin{equation}
 \cF(\widetilde \bw_{R,R'};B_{R'})
 =\cF(\bw_R;B_R)+\cp(R^{2-p}-R'^{2-p}).
\end{equation}
Using this competitor and the minimality of $\bw_R$ gives
\begin{equation}\label{eq:mono}
 \J(R')
=\min_{\bv\in\cA_{R'}}\cF(\bv;B_{R'})+\cp R'^{2-p} \leqslant\min_{\bv\in\cA_R}\cF(\bv;B_R)+\cp R^{2-p}
 =\J(R).
\end{equation}
Consequently, $\J$ is decreasing and, being nonnegative,
has a finite limit as $R\to\infty$.

We next identify $\lim_{R\to\infty}\J(R)$ with the whole-plane
energy infimum at degree one. We first define degree at
infinity and the corresponding admissible classes.
For every $\bv\in W^{1,p}_{\mathrm{loc}}(\R^2;\R^2)$ with
$\cF(\bv;\R^2)<\infty$, Morrey's inequality and the
integrability of the potential imply that
$|\bv(y)|\to1$ uniformly as $|y|\to\infty$.\footnote{Indeed,
the global Morrey estimate gives
$|\bv(x)-\bv(y)|\leqslant
C_p\|\nabla \bv\|_{L^p(\R^2)}|x-y|^{1-2/p}$.
If the asserted convergence failed, there would exist
$\eta>0$ and points $y_k$ with $|y_k|\to\infty$ such that
$\bigl||\bv(y_k)|-1\bigr|\geqslant\eta$.
The Morrey estimate then gives a radius $s>0$,
independent of $k$, on which
$\bigl||\bv|-1\bigr|\geqslant\eta/2$ throughout $B_s(y_k)$.
After extracting disjoint balls, each contributes at
least $\pi s^2\eta^2/16>0$ to the potential energy,
contradicting $\cF(\bv;\R^2)<\infty$.
See also \cite[Lemma~2.1]{ABGS}.}
Consequently, $\bv$ has no zeros outside some disk.
For any two sufficiently large circles centered at the origin,
varying the radius gives a homotopy between the restrictions of
$\bv/|\bv|$. Their degrees therefore agree, and we denote
this common integer by $\deg_\infty \bv$.
Using this degree, for each $d\in\mathbb Z$ define
\begin{equation}\label{eq:wp}
 \cX_d=
 \left\{\bv\in W^{1,p}_{\mathrm{loc}}(\R^2;\R^2):
 \cF(\bv;\R^2)<\infty,\ \deg_\infty \bv=d\right\},
 \qquad
\mdeg_d=\inf_{\bv\in\cX_d}\cF(\bv;\R^2).
\end{equation}
Because $\cX_1$ relaxes the exterior radial constraint to a degree-$1$ condition at $\infty$, $\mdeg_1$ is the natural candidate for $\lim_{R\to\infty} \J(R)$. The following lemma identifies this limit.

\begin{lemma}[Identification of the limiting energy]
\label{lem:core}
For every fixed $p>2$, the whole-plane infimum $\mdeg_1$, defined by \eqref{eq:wp} with $d=1$, satisfies
\begin{equation}\label{eq:core}
0<\mdeg_1=\lim_{R\to\infty}\J(R)<\infty.
\end{equation}
Thus the prescribed radial trace and the degree-one condition at
infinity give the same limiting infimum.
\end{lemma}

\begin{proof}[Proof of Lemma~\ref{lem:core}]
We approximate an arbitrary finite-energy degree-one map
by maps with the prescribed radial trace on large circles.
We first restore the modulus to one in an annulus of
width one, then remove the phase correction in a larger
annulus. Finally, we prove a uniform positive lower
bound for the energy of a degree-one map.

Fix $p\in(2,\infty)$. For every $R>0$, extending a
minimizer $\bw_R\in\cA_R$ by $\bn$
gives a map in $\cX_1$ with energy
$\J(R)$. Consequently,
\begin{equation}\label{eq:core_lb}
 0\leqslant\mdeg_1
 \leqslant\J(R)<\infty,
 \qquad R>0.
\end{equation}
Finiteness follows, for example, by using
$x\mapsto x/R\in W^{1,p}(B_R;\R^2)$ as a competitor
in $\cA_R$. By \eqref{eq:mono},
$\J$ is nonincreasing, so
$\lim_{R\to\infty}\J(R)$ exists
and is finite.

Fix $\bv\in\cX_1
\subset W^{1,p}_{\mathrm{loc}}(\R^2;\R^2)$.
The uniform convergence $|\bv|\to1$ at infinity,
established before \eqref{eq:wp},
allows us to choose radii $R_j\in[1,\infty)$,
$j\in\mathbb N$, such that
\begin{equation}\label{eq:radii}
 R_j\longrightarrow\infty,
 \qquad
 \frac12\leqslant|\bv|\leqslant\frac32
 \quad\text{on }\R^2\setminus B_{R_j},
 \qquad
 \deg\left(\frac{\bv}{|\bv|},\partial B_r\right)=1
 \quad(r\geqslant R_j).
\end{equation}
Define the tail quantities
\begin{equation}\label{eq:core_tail}
 t_j:=
 \int_{\R^2\setminus B_{R_j}}
 \bigl(|\nabla\bv|^p+(1-|\bv|^2)^2\bigr)\dd
 \in[0,\infty),
 \qquad
 t_j\longrightarrow0
 \quad\text{as }j\to\infty.
\end{equation}
Polar integration and Sobolev slicing give
$r_j\in(R_j,2R_j)$ such that
\begin{equation}\label{eq:core_trace}
 \bh_j\in W^{1,p}(\TT;\C),
 \qquad
 \bh_j(\theta):=\bv(r_je^{i\theta}),
 \qquad \theta\in\TT,
\end{equation}
and
\begin{equation}\label{eq:circle_bound}
 \int_{\partial B_{r_j}}
 \bigl(|\partial_\tau\bv|^p
       +(1-|\bv|^2)^2\bigr)\,d\mathcal H^1
 \leqslant\frac{2t_j}{R_j}.
\end{equation}
Here $\partial_\tau$ denotes differentiation along the
positively oriented unit tangent to the circle.
All angular integrals below are over $[0,2\pi]$,
with matching endpoint traces. Set
\begin{equation}\label{eq:norm_trace}
 \rh_j:=|\bh_j|
 \in W^{1,p}(\TT;\R),
 \qquad
 \bm g_j:=\frac{\bh_j}{\rh_j}
 \in W^{1,p}(\TT;\Sone).
\end{equation}
The Sobolev chain rule and
\eqref{eq:radii} imply
\begin{equation}\label{eq:trace_deriv}
 |\rh_j'|\leqslant|\bh_j'|,
 \qquad
 |\bm g_j'|\leqslant
 \rh_j^{-1}|\bh_j'|
 \overset{\eqref{eq:radii}}{\leqslant}
 2|\bh_j'|
 \quad\text{a.e. on }\TT.
\end{equation}
We first restore the modulus to one. Define
$\bm A_j\in
W^{1,p}(B_{r_j+1}\setminus\overline{B_{r_j}};\C)$ by
\begin{equation}\label{eq:mod_fill}
 \bm A_j(re^{i\theta})
 :=\bigl((1-s)\rh_j(\theta)+s\bigr)\bm g_j(\theta),
 \qquad
 s=r-r_j,
 \qquad
 (r,\theta)\in(r_j,r_j+1)\times\TT.
\end{equation}
Its continuous representative has boundary values
\begin{equation}\label{eq:mod_trace}
 (\bm A_j)_{|\partial B_{r_j}}(r_je^{i\theta})
 =\bh_j(\theta),
 \qquad
 (\bm A_j)_{|\partial B_{r_j+1}}((r_j+1)e^{i\theta})
 =\bm g_j(\theta),
 \qquad \theta\in\TT.
\end{equation}
Differentiating \eqref{eq:mod_fill}
and using \eqref{eq:trace_deriv} gives
\begin{equation}\label{eq:mod_deriv}
 |\partial_r\bm A_j|=|1-\rh_j|,
 \qquad
 |\partial_\theta\bm A_j|
 \leqslant|\rh_j'|+\frac32|\bm g_j'|
 \overset{\eqref{eq:trace_deriv}}{\leqslant}
 4|\bh_j'|
\end{equation}
almost everywhere in polar coordinates.
Since $|\bm A_j|$ lies between $\rh_j$ and one,
and $|1-\rh_j|\leqslant1/2$, we also have
\begin{equation}\label{eq:mod_pot}
 (1-|\bm A_j|^2)^2\leqslant(1-\rh_j^2)^2,
 \qquad
 |1-\rh_j|^p
 \leqslant|1-\rh_j|^2
 \leqslant(1-\rh_j^2)^2.
\end{equation}
Using \eqref{eq:mod_deriv},
\eqref{eq:mod_pot}, and
$r_j\leqslant r\leqslant r_j+1\leqslant2r_j$
in polar coordinates yields
\begin{equation}\label{eq:mod_en}
 \cF(\bm A_j;
 B_{r_j+1}\setminus\overline{B_{r_j}})
 \leqslant C_p\int_{\partial B_{r_j}}
 \bigl(|\partial_\tau\bv|^p
       +(1-|\bv|^2)^2\bigr)\,d\mathcal H^1
 \overset{\eqref{eq:circle_bound}}{\leqslant}
 C_p\frac{t_j}{R_j},
\end{equation}
where $C_p\in(0,\infty)$ depends only on $p$.

We next remove the phase correction.
The continuous map $\bm g_j$ has degree one.
Lifting its phase on $[0,2\pi]$ gives
\begin{equation}\label{eq:phase_lift}
 \bm g_j(\theta)
 =e^{i(\theta+\ps_j(\theta))},
 \qquad
 \ps_j\in W^{1,p}(\TT;\R),
 \qquad
 |\langle\ps_j\rangle|\leqslant\pi.
\end{equation}
Indeed, integrating $\bm g_j\wedge\bm g_j'$ and
choosing the initial phase to match $\bm g_j(0)$
constructs a lift in $W^{1,p}(0,2\pi;\R)$.
Its increment is $2\pi$, so subtracting $\theta$
gives a periodic phase correction.
Adding a multiple of $2\pi$ then places its mean
in $[-\pi,\pi]$.
The one-dimensional Poincar\'e inequality and
\eqref{eq:trace_deriv} imply
\begin{equation}\label{eq:phase_bound}
 \int_0^{2\pi}
 \bigl(|\ps_j|^p+|\ps_j'|^p\bigr)
 \,d\theta
 \leqslant C_p\left(
 1+\int_0^{2\pi}|\bh_j'|^p\,d\theta\right).
\end{equation}
Set $b_j:=r_j+1\in(1,\infty)$ and define
$\bm B_j\in
W^{1,p}(B_{2b_j}\setminus\overline{B_{b_j}};\Sone)$
by
\begin{equation}\label{eq:phase_interp}
 \bm B_j(re^{i\theta})
 :=\exp\left(
 i\left[\theta+
 \left(2-\frac r{b_j}\right)\ps_j(\theta)\right]
 \right),
 \qquad
 (r,\theta)\in(b_j,2b_j)\times\TT.
\end{equation}
Its continuous representative satisfies
\begin{equation}\label{eq:phase_trace}
 (\bm B_j)_{|\partial B_{b_j}}(b_je^{i\theta})
 =\bm g_j(\theta),
 \qquad
 (\bm B_j)_{|\partial B_{2b_j}}(2b_je^{i\theta})
 =e^{i\theta}=\bn(2b_je^{i\theta}),
 \qquad \theta\in\TT.
\end{equation}
The map is unit valued, so its potential energy vanishes.
Its radial and angular phase derivatives are
$-\ps_j/b_j$ and
$1+(2-r/b_j)\ps_j'$, respectively.
Moreover, the tangential chain rule gives
\begin{equation}\label{eq:ang_en}
 \int_0^{2\pi}|\bh_j'|^p\,d\theta
 =r_j^{p-1}
 \int_{\partial B_{r_j}}
 |\partial_\tau\bv|^p\,d\mathcal H^1.
\end{equation}
Polar integration, followed by
\eqref{eq:phase_bound},
\eqref{eq:ang_en}, and
\eqref{eq:circle_bound}, therefore yields
\begin{equation}\label{eq:phase_ann}
 \cF(\bm B_j;
 B_{2b_j}\setminus\overline{B_{b_j}})
 \leqslant C_pb_j^{2-p}
 \left(1+\int_0^{2\pi}|\bh_j'|^p\,d\theta\right)
 \leqslant C_p(R_j^{2-p}+t_j).
\end{equation}
For the last inequality, we used
$b_j\geqslant r_j\geqslant R_j$,
$r_j\leqslant2R_j$, and $p>2$.

Set $\widehat R_j:=2b_j$, so that
$\widehat R_j\to\infty$, and glue the three maps:
\begin{equation}\label{eq:glue}
 \bv_j(x):=
 \begin{cases}
 \bv(x),&x\in B_{r_j},\\
 \bm A_j(x),&x\in B_{b_j}\setminus\overline{B_{r_j}},\\
 \bm B_j(x),&x\in B_{\widehat R_j}\setminus\overline{B_{b_j}}.
 \end{cases}
\end{equation}
The trace identities
\eqref{eq:mod_trace} and
\eqref{eq:phase_trace} give
$\bv_j\in\cA_{\widehat R_j}
\subset W^{1,p}(B_{\widehat R_j};\R^2)$
by Sobolev gluing.
With a sufficiently large constant $C_p>0$, define
\begin{equation}\label{eq:trunc_err}
 \delta_j:=
 C_p\left(\frac{t_j}{R_j}+t_j+R_j^{2-p}\right)>0,
 \qquad
 \delta_j\longrightarrow0
 \quad\text{as }j\to\infty.
\end{equation}
The estimates
\eqref{eq:mod_en} and
\eqref{eq:phase_ann} now give
\begin{equation}\label{eq:trunc}
 \cF(\bv_j;B_{\widehat R_j})
 \leqslant\cF(\bv;\R^2)+\delta_j,
 \qquad
 \widehat R_j\longrightarrow\infty
 \quad\text{as }j\to\infty.
\end{equation}
By the definition \eqref{eq:j} of $\J$,
\begin{equation}\label{eq:core_ub}
 \begin{split}
 \J(\widehat R_j)
 &\overset{\eqref{eq:j}}{\leqslant}
 \cF(\bv_j;B_{\widehat R_j})
 +\cp\widehat R_j^{2-p}\\
 &\overset{\eqref{eq:trunc}}{\leqslant}
 \cF(\bv;\R^2)+\delta_j
 +\cp\widehat R_j^{2-p}.
 \end{split}
\end{equation}
Letting $j\to\infty$ and then taking the infimum over
$\bv\in\cX_1$ yields
\begin{equation}\label{eq:core_lim}
 \lim_{R\to\infty}\J(R)
 \leqslant\mdeg_1.
\end{equation}
Together with \eqref{eq:core_lb},
this identifies the limit.

It remains to prove strict positivity.
Every $\bv\in\cX_1$ has a zero
$a\in\R^2$: otherwise its normalization would
extend a degree-one boundary map continuously
across a sufficiently large disk.
Write $M_v:=\cF(\bv;\R^2)\in[0,\infty)$.
Morrey's inequality on $B_1(a)$, applied to differences
from $\bv(a)=0$, gives a constant
$C_M(p)\in(0,\infty)$, depending only on $p$, such that
\begin{equation}\label{eq:core_morrey}
 \sup_{B_1(a)}|\bv|
 \leqslant C_M(p)
 \|\nabla\bv\|_{L^p(B_1(a);\R^{2\times2})}
 \leqslant C_M(p)(pM_v)^{1/p}.
\end{equation}
Set
\begin{equation}\label{eq:core_small}
 \mu_p:=\frac1{p(2C_M(p))^p}\in(0,\infty).
\end{equation}
If $M_v\leqslant\mu_p$, then
\eqref{eq:core_morrey} gives
$|\bv|\leqslant1/2$ on $B_1(a)$, and hence
\begin{equation}\label{eq:pot_lb}
 M_v\geqslant
 \frac14\int_{B_1(a)}(1-|\bv|^2)^2\dd
 \geqslant\frac{9\pi}{64}.
\end{equation}
If $M_v>\mu_p$, this already supplies a positive
lower bound. Taking the infimum over
$\bv\in\cX_1$ therefore gives
\begin{equation}\label{eq:core_pos}
 \mdeg_1
 \geqslant\min\left\{\mu_p,\frac{9\pi}{64}\right\}>0.
\end{equation}
This proves the strict positivity and completes
\eqref{eq:core}.
\end{proof}

We next use the energy limit in
Lemma~\ref{lem:core} to control the region
where the modulus stays a fixed amount below one.
For $\eta\in(0,1]$ and $\e>0$, we call
\begin{equation}
 D_{\eta}(\bu_\e)
 :=\{x\in B_1:|\bu_\e(x)|\leqslant1-\eta\}
\end{equation}
the \emph{modulus-defect set at threshold $\eta$}.

For each fixed $\eta\in(0,1)$, the energy identification
and compactness up to translations will be used in
the proof of Proposition~\ref{prop:conc} to exclude
families of global minimizers $\bu_\e\in\cA$
and points $x_\e,y_\e\in D_{\eta}(\bu_\e)$ such that
\begin{equation}
 \frac{|x_\e-y_\e|}{\lc}
 \longrightarrow\infty,
 \qquad \lc=\e^{2/p},
 \qquad \e\downarrow0.
\end{equation}
Geometrically, this limit would mean that two points
with the same fixed loss of modulus move arbitrarily
far apart in the rescaled coordinates $y=x/\lc$,
even if their distance in the original disk tends
to zero. Excluding this behavior prevents the defect
set from splitting into regions separated by distances
much larger than the natural scale $\lc$.

Proposition~\ref{prop:conc} makes this exclusion
quantitative: the entire set $D_{\eta}(\bu_\e)$ remains
within a bounded distance of any chosen zero in the
rescaled coordinates. The bound and the smallness
condition on $\e$ depend only on $p$ and $\eta$,
uniformly over all global minimizers and all choices
of the zero.

\begin{proposition}[Concentration at the natural scale]
\label{prop:conc}
For every $p>2$ and $\eta\in(0,1]$, there exist constants
$C_{p,\eta}>0$ and $\e_{p,\eta}\in(0,1)$ such that,
for every $\e\in(0,\e_{p,\eta})$, every global minimizer
$\bu_\e\in\cA$ of $\E_\e$, and every $a_\e\in B_1$
satisfying $\bu_\e(a_\e)=0$, one has
\begin{equation}\label{eq:conc}
D_{\eta}(\bu_\e):= \{x\in B_1:|\bu_\e(x)|\leqslant1-\eta\}
 \subset B_{C_{p,\eta}\lc}(a_\e),
 \qquad \lc=\e^{2/p}.
\end{equation}
The constants $C_{p,\eta}$ and $\e_{p,\eta}$ depend
only on $p$ and $\eta$; in particular, they are independent
of the chosen minimizer $\bu_\e$ and its zero $a_\e$.
Thus the concentration estimate is uniform over all
global minimizers as $\e\downarrow0$: the modulus-defect set has diameter at most
$2C_{p,\eta}\e^{2/p}$, which tends to zero as $\e\downarrow0$.
\end{proposition}

\begin{remark}[The threshold endpoints]
The condition $\eta>0$ is essential. At $\eta=0$,
the bound $|\bu_\e|\leqslant1$ makes the modulus-defect
set equal to $B_1$, so concentration in a ball of
radius $O(\lc)$ is impossible as $\e\downarrow0$.
The endpoint $\eta=1$, however, is admissible:
the modulus-defect set then equals the zero set
$\Zset_\e=\{x\in B_1:\bu_\e(x)=0\}$.
Since $\Zset_\e\subset\{|\bu_\e|\leqslant1/2\}$,
Proposition~\ref{prop:conc} gives
\begin{equation}
 \Zset_\e\subset B_{C_{p,1/2}\e^{2/p}}(a_\e)
 \qquad\text{for every }a_\e\in \Zset_\e,
 \qquad 0<\e<\e_{p,1/2}.
\end{equation}
In particular, for every prescribed radius $s>0$,
there exists $\e_{p,s}>0$ such that
$\Zset_\e\subset B_s(a_\e)$ for every
$\e\in(0,\e_{p,s})$, every global minimizer $\bu_\e$,
and every choice of its zero $a_\e$.
Thus any zero can serve as the center of an arbitrarily
small ball containing all the others, provided $\e$
is sufficiently small. This expresses uniform
concentration of the zero set as $\e\downarrow0$,
but does not imply uniqueness or isolation of its points
for a fixed $\e$.
Its concentration follows directly from the estimate
at $\eta=1/2$, with the same constants and smallness
restriction on $\e$.
\end{remark}

\begin{remark}
The containing ball may intersect $\partial B_1$:
Proposition~\ref{prop:conc} does not yet locate
$a_\e$, whose convergence to the origin will be proved
in Section~\ref{sec:proof}.
Outside this ball,
$|\bu_\e|>1-\eta$.
This lower bound permits normalization of $\bu_\e$,
supplies the modulus condition for the filling lemma
in Section~\ref{sec:fill}, and allows the phase
comparison in Section~\ref{sec:defect}.
\end{remark}

\begin{proof}[Proof of Proposition~\ref{prop:conc}]
We first prove compactness up to translations for
whole-plane minimizing sequences. A nonzero degree
costs at least $\mdeg_1$, leaving no energy
for another separated modulus defect. We then apply
this compactness to the rescaled disk minimizers.

Fix $p\in(2,\infty)$. We use the degree-reduction inequality
\begin{equation}\label{eq:deg_red}
 \cF(\bv;\R^2)
 \geqslant\mdeg_1,
 \qquad
 d\in\mathbb Z\setminus\{0\},
 \qquad
 \bv\in\cX_d.
\end{equation}
For $d\geqslant2$, \cite[Proposition~1]{ABGS}
constructs a degree-one map with the same first component
and the same absolute value of the second component.
It therefore preserves $|\bv|$ and
$|\nabla\bv|$ almost everywhere, so its construction
applies to our normalization. Negative degrees are
handled by complex conjugation; degree one is covered
by the definition of $\mdeg_1$.

\smallskip\noindent
\emph{Step 1: compactness up to translations.}
Let $(\bV_k)_{k\in\mathbb N}
\subset\cX_1$ satisfy
\begin{equation}\label{eq:min_seq}
 |\bV_k|\leqslant1,
 \qquad
 \cF(\bV_k;\R^2)
 \longrightarrow\mdeg_1
 \quad\text{as }k\to\infty.
\end{equation}
After discarding finitely many terms, the energies
are bounded by $\mdeg_1+1$.
The global Morrey inequality gives $s\in(0,\infty)$,
depending only on $p$, such that
\begin{equation}\label{eq:osc}
 |\bV_k(x)-\bV_k(y)|
 \leqslant\frac14
 \qquad
 (x,y\in\R^2,\ |x-y|\leqslant s,\ k\in\mathbb N).
\end{equation}
Thus every point of the compact nonempty set
$S_k:=\{x\in\R^2:|\bV_k(x)|\leqslant1/2\}$
has a surrounding disk on which the modulus is at most
$3/4$. In particular,
\begin{equation}\label{eq:loc_cost}
 \cF(\bV_k;B_s(a))
 \geqslant
 \frac{\pi s^2}{4}\bigl[1-(3/4)^2\bigr]^2
 =:\gamma_0>0,
 \qquad a\in S_k.
\end{equation}
Nonemptiness follows from degree one; compactness
follows from $|\bV_k|\to1$ at infinity.

Choose centers in $S_k$ successively, each outside
the closed $3s$-disks around the preceding centers.
Their open $s$-disks are disjoint, so
\eqref{eq:loc_cost} bounds their number
uniformly. Passing to a subsequence, we obtain
$N\in\mathbb N$ centers
$x_{1,k},\ldots,x_{N,k}\in S_k$ whose closed
$3s$-disks cover $S_k$.

Pass to a further subsequence on which every pairwise
distance has a limit in $[0,\infty]$. Group indices
whose centers remain at bounded distance; the triangle
inequality makes this an equivalence relation.
Let $K\in\{1,\ldots,N\}$ be the number of groups and
$c_{j,k}\in S_k$ a representative of group $j$.
There exists $L\in(0,\infty)$ such that
\begin{equation}\label{eq:cover}
 S_k\subset\bigcup_{j=1}^K B_L(c_{j,k}),
 \qquad
 |c_{j,k}-c_{h,k}|\longrightarrow\infty
 \quad(j\ne h)
 \quad\text{as }k\to\infty.
\end{equation}
Here $L$ includes the covering radii and all bounded
distances within each group.

Local weak compactness and the common Morrey bound
give maps
$\bm v_j\in W^{1,p}_{\mathrm{loc}}(\R^2;\R^2)$ such that,
along a common subsequence,
\begin{equation}\label{eq:prof_conv}
 \bV_k(c_{j,k}+\cdot)
 \rightharpoonup\bm v_j
 \quad\text{in }W^{1,p}_{\mathrm{loc}}(\R^2;\R^2),
 \qquad j\in\{1,\ldots,K\},
\end{equation}
with uniform convergence on compact sets as well.
Lower semicontinuity gives finite energy.
Moreover, \eqref{eq:osc}
and $c_{j,k}\in S_k$ imply $|\bm v_j|\leqslant3/4$
on $B_s$, hence
$\cF(\bm v_j;\R^2)\geqslant\gamma_0$.

For each fixed $T>L$, the circles
$\partial B_T(c_{j,k})$ lie outside all covering balls
for large $k$, so $|\bV_k|>1/2$ there.
The same observation at every fixed point outside
$B_L$ gives $|\bm v_j|\geqslant1/2$ on
$\R^2\setminus B_L$. Thus
$\bm v_j\in\cX_{d_j}$ for some
$d_j\in\mathbb Z$.
Uniform convergence on $\partial B_T$ preserves the
degrees of the normalized maps. Since the disjoint
disks $B_T(c_{j,k})$ contain every zero, degree
additivity gives
\begin{equation}\label{eq:deg_sum}
 \sum_{j=1}^K d_j=\deg_\infty\bV_k=1
 \qquad\text{for all sufficiently large }k.
\end{equation}
Indeed, the normalization is continuous on the
complement of these disks inside any sufficiently
large disk; its outer degree equals the sum of
the inner degrees, with all circles positively oriented.

For fixed $T>L$, the disks $B_T(c_{j,k})$ are disjoint
for large $k$. Lower semicontinuity, followed by
$T\to\infty$, therefore yields
\begin{equation}\label{eq:prof_en}
 \sum_{j=1}^K\cF(\bm v_j;\R^2)
 \leqslant
 \liminf_{k\to\infty}\cF(\bV_k;\R^2)
 \overset{\eqref{eq:min_seq}}{=}
 \mdeg_1.
\end{equation}
At least one $d_j$ is nonzero by
\eqref{eq:deg_sum}.
Its energy is at least $\mdeg_1$
by \eqref{eq:deg_red};
each remaining profile costs at least $\gamma_0$.
Consequently,
\begin{equation}\label{eq:one_prof}
 \mdeg_1+(K-1)\gamma_0
 \leqslant
 \sum_{j=1}^K\cF(\bm v_j;\R^2)
 \overset{\eqref{eq:prof_en}}{\leqslant}
 \mdeg_1.
\end{equation}
Thus $K=1$, $d_1=1$, and
$\cF(\bm v_1;\R^2)=\mdeg_1$.
Writing $c_k:=c_{1,k}$ and $\bv:=\bm v_1$, we have proved
\begin{equation}\label{eq:compact}
 \bV_k(c_k+\cdot)\rightharpoonup\bv
 \quad\text{in }W^{1,p}_{\mathrm{loc}}(\R^2;\R^2),
 \qquad
 \bv\in\cX_1,
 \qquad
 \cF(\bv;\R^2)=\mdeg_1,
\end{equation}
with local uniform convergence as $k\to\infty$.

\smallskip\noindent
\emph{Step 2: uniform concentration of the defect sets.}
Fix $\eta\in(0,1)$.
For an open set $G\subset\R^2$ and
$\bm W\in W^{1,p}_{\mathrm{loc}}(G;\R^2)$, extend
the defect-set notation by
\begin{equation}\label{eq:def_set}
 \Ddef_\eta(\bm W)
 :=\{y\in G:|\bm W(y)|\leqslant1-\eta\}.
\end{equation}
Suppose that the rescaled defect diameters of disk
minimizers admit no uniform bound as $\e\downarrow0$.
For each $k\in\mathbb N$, choose $\e_k\in(0,1/k)$
and a global minimizer
$\bu_{\e_k}\in\cA$ whose rescaled defect diameter
exceeds $k$. Set $\lc_k:=\e_k^{2/p}$ and extend
its rescaling by $\bn$:
\begin{equation}\label{eq:rescaled}
 \bW_k(y):=
 \begin{cases}
 \bu_{\e_k}(\lc_k y),
     &y\in B_{\lc_k^{-1}},\\
 \bn(y),
     &y\in\R^2\setminus B_{\lc_k^{-1}}.
 \end{cases}
\end{equation}
Matching traces give
$\bW_k\in\cX_1
\subset W^{1,p}_{\mathrm{loc}}(\R^2;\R^2)$,
with $|\bW_k|\leqslant1$.
Lemma~\ref{lem:core} and our choice of minimizers give
\begin{equation}\label{eq:rescaled_en}
 \cF(\bW_k;\R^2)
 =\J(\lc_k^{-1})
 \longrightarrow\mdeg_1,
 \qquad
 \operatorname{diam}\Ddef_\eta(\bW_k)
 \longrightarrow\infty
 \quad\text{as }k\to\infty.
\end{equation}
Apply \eqref{eq:compact} and set
$\widehat{\bW}_k
:=\bW_k(c_k+\cdot)\in\cX_1$.
Translation preserves diameter, so we can choose
$b_k\in\Ddef_\eta(\widehat{\bW}_k)
\subset\R^2$ with $|b_k|\to\infty$.
Indeed, a bounded set of diameter $d$ contains a point
whose distance from the origin is at least $d/2$.

The uniform energy bound gives a Morrey constant
$\Hmor\in(0,\infty)$, depending only on $p$,
such that
\begin{equation}\label{eq:morrey}
 |\widehat{\bW}_k(x)
   -\widehat{\bW}_k(y)|
 \leqslant\Hmor|x-y|^{1-2/p},
 \qquad x,y\in\R^2,\quad k\in\mathbb N.
\end{equation}
Choose $\rpack\in(0,\infty)$ with
$\Hmor\rpack^{1-2/p}\leqslant\eta/2$.
Then $|\widehat{\bW}_k|\leqslant1-\eta/2$
on $B_{\rpack}(b_k)$, and
\begin{equation}\label{eq:def_cost}
 \cF(\widehat{\bW}_k;
 B_{\rpack}(b_k))
 \geqslant\frac{\pi\rpack^2}{4}
 \bigl[1-(1-\eta/2)^2\bigr]^2
 =:\gamma>0.
\end{equation}
For every fixed $T>0$, this disk is disjoint from $B_T$
for large $k$. Nonnegativity of the energy density,
\eqref{eq:def_cost}, and lower semicontinuity give
\begin{equation}\label{eq:escape_en}
 \mdeg_1
 \overset{\eqref{eq:rescaled_en}}{=}
 \lim_{k\to\infty}
 \cF(\widehat{\bW}_k;\R^2)
 \geqslant\cF(\bv;B_T)+\gamma.
\end{equation}
Letting $T\to\infty$ contradicts
\eqref{eq:compact}:
\begin{equation}\label{eq:no_split}
 \mdeg_1
 \geqslant\cF(\bv;\R^2)+\gamma
 \overset{\eqref{eq:compact}}{=}
 \mdeg_1+\gamma>\mdeg_1.
\end{equation}

Thus a uniform diameter bound and a smallness threshold
exist, depending only on $p$ and $\eta$.
Every zero belongs to the defect set.
Taking any zero as center, rescaling back, and choosing
the radius constant strictly larger than the diameter
bound gives \eqref{eq:conc}.
For $\eta=1$, apply the result at $\eta=1/2$, since the
zero set is contained in that defect set.
All choices are uniform over the minimizers and their zeros.
\end{proof}

\begin{remark}[Dependence of the concentration radius]
The proof of Proposition~\ref{prop:conc} gives no
explicit formula for $C_{p,\eta}$. Morrey's inequality
and the energy bound quantitatively control the number
of disjoint balls centered at modulus-defect points.
Their mutual distances, however, are bounded only
through the compactness contradiction in
\eqref{eq:no_split}, which yields no explicit estimate.
For fixed $p$, increasing $\eta$ shrinks the set
$\{|\bu_\e|\leqslant1-\eta\}$. Hence any constant and
smallness restriction valid at $\eta_0\in(0,1)$ remain
valid for every $\eta\in[\eta_0,1)$.
But the argument provides neither a uniform bound as
$\eta\downarrow0$ nor monotonicity in $p$.
\end{remark}

\section{Improving a nonradial exterior trace}\label{sec:fill}

To obtain the exterior energy bound in
Section~\ref{sec:excess}, we construct a competitor
by centering the map at a chosen zero $a_\e$ and
interpolating between its boundary trace on $\partial B_\delta$,
$0<\delta<1$, and $\bn$ on $\partial B_1$.
The comparison requires choosing a radius $\delta>0$
that satisfies two competing constraints: $\delta$
must be small enough that the domain of the
translated map $x\mapsto\bu_\e(a_\e+x)$ is
contained in $B_1$, yet large enough that
$|\bu_\e|$ is bounded away from zero on
$\partial B_\delta(a_\e)$, allowing us to normalize
its trace on $\partial B_\delta$.

More precisely, fix $\e\in(0,\e_{p,1/2})$, a global
minimizer $\bu_\e\in\cA$, a zero
$a_\e\in B_1$, and set $\lc=\e^{2/p}$.
Applying Proposition~\ref{prop:conc} with
threshold $\eta=1/2$ yields
\begin{equation}
D_{1/2}(\bu_\e)
=\{x\in B_1:|\bu_\e(x)|\leqslant1/2\}
\subset B_{C_{p,1/2}\lc}(a_\e).
\end{equation}
Consequently, whenever $\delta>C_{p,1/2}\lc$,
we have
\begin{equation}
	|\bu_\e(x)|>\frac12
	\qquad\text{for all }
	x\in\partial B_\delta(a_\e)\cap B_1.
\end{equation}
A priori, this supplies a nonvanishing trace only on
the portion of the comparison circle lying inside
$B_1$. If $C_{p,1/2}\lc<\delta<1-|a_\e|$,
then $\overline{B_\delta(a_\e)}\subset B_1$, ensuring
that the trace
$(\bu_\e)_{|\partial B_\delta(a_\e)}$
is well-defined on the entire circle and has modulus
strictly greater than $1/2$.
At this stage, however, no estimate guarantees that
the interval $(C_{p,1/2}\lc,1-|a_\e|)$ is
nonempty; thus, for a given
$\delta>C_{p,1/2}\lc$, the translated map
$x\mapsto\bu_\e(a_\e+x)$ need not be defined on
all of $B_\delta$. To circumvent this issue, we
extend $\bu_\e$ to all of $\R^2$ prior to
selecting the comparison circle:
\begin{equation}\label{eq:ext}
	\bU_\e\in W^{1,p}_{\mathrm{loc}}(\R^2;\R^2),
	\qquad
	\bU_\e(x)=
	\begin{cases}
		\bu_\e(x),&x\in B_1,\\
		\bn(x),&x\in\R^2\setminus B_1.
	\end{cases}
\end{equation}
Since the traces match on $\partial B_1$, Sobolev gluing yields
$\bU_\e\in W^{1,p}_{\mathrm{loc}}(\R^2;\R^2)$.
For $p>2$, we may pass to the continuous representative of $\bU_\e$;
its normalization $\bU_\e/|\bU_\e|$
is then continuous away from the zero set, which allows
us to compare the topological degrees on surrounding
circles via homotopy invariance.

Take $\e<\e_{p,1/2}$ small enough that $C_{p,1/2}\lc<1$.
Combining the concentration estimate of Proposition~\ref{prop:conc}
inside $B_1$ with the identity $|\bU_\e|=|\bn|=1$
outside $B_1$ from \eqref{eq:ext}, we obtain the
global pointwise bound
\begin{equation}
	\frac12\leqslant|\bU_\e(x)|\leqslant1
	\qquad\text{for all }
	x\in\R^2\setminus B_{C_{p,1/2}\lc}(a_\e).
\end{equation}
In particular, for every $\delta>C_{p,1/2}\lc$, the translated
map $x\mapsto\bU_\e(a_\e+x)$ is well-defined on all of $B_\delta$
and $\bU_\e$ is nonvanishing on $\partial B_\delta(a_\e)$,
regardless of whether $\partial B_\delta(a_\e)$ is contained in $B_1$.
This allows us to perform the exterior comparison
without assuming a priori that $a_\e\to 0$, a fact
we will establish later in Section~\ref{sec:proof}.

For $\delta\in(C_{p,1/2}\lc,1)$, we define the parametrized boundary trace
\begin{equation}
	\bh_{\e,\delta}(\omega)
	:=\bU_\e(a_\e+\delta\omega),
	\qquad \omega\in\Sone.
\end{equation}
By Sobolev embedding and the trace theorem,
$\bh_{\e,\delta}\in C^0(\Sone;\R^2)\cap W^{1-1/p,p}(\Sone;\R^2)$
for every such $\delta$, while Fubini's theorem
(or standard Sobolev slicing) yields the stronger regularity
\begin{equation}
	\bh_{\e,\delta}\in W^{1,p}(\Sone;\R^2)
	\qquad\text{for almost every }
	\delta\in(C_{p,1/2}\lc,1).
\end{equation}

In what follows, we will work with such radii\footnote{The restriction to almost every radius is needed
to ensure that $\bh_{\e,\delta}\in W^{1,p}(\Sone;\R^2)$. This additional regularity makes
the angular energy finite and allows us to apply
Lemma~\ref{lem:fill}, as continuity alone does not
control the tangential derivative of the trace.}  $\delta$.
For every $\delta\in(C_{p,1/2}\lc,1)$, the
continuous trace satisfies
\begin{equation}
	\frac12\leqslant|\bh_{\e,\delta}|\leqslant1
	\quad\text{on }\Sone,
	\qquad
	\deg\left(
	\frac{\bh_{\e,\delta}}{|\bh_{\e,\delta}|},
	\Sone\right)=1.
\end{equation}
Indeed, since $\bU_\e/|\bU_\e|$ is continuous
on $\R^2\setminus B_\delta(a_\e)$, homotopy invariance implies
that its degree on $\partial B_\delta(a_\e)$ agrees with its degree
on any large circle enclosing $B_1$, where $\bU_\e$ coincides
with $\bn$ (with all circles positively oriented).

\begin{figure}[tbp]
	\centering
	\resizebox{\linewidth}{!}{%
		\includegraphics[]{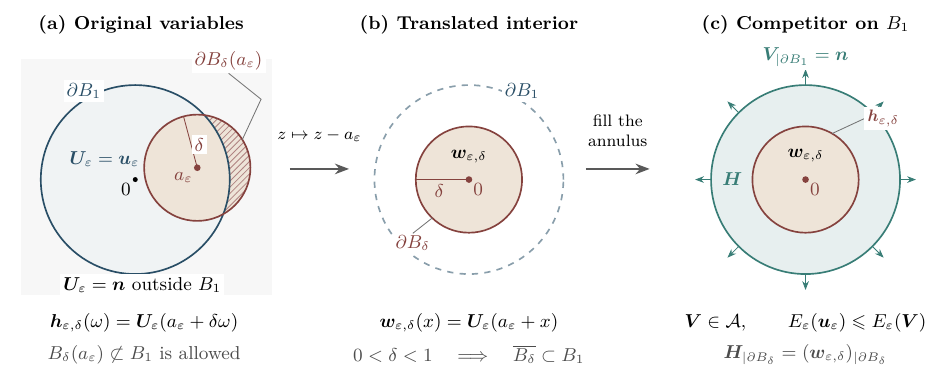}%
	}
	\caption{Translation and filling. The extension
		$\bU_\e=\bn$ supplies the values
		on the hatched part of $B_\delta(a_\e)$ outside
		$B_1$. Translation gives $\bw_{\e,\delta}$
		on $B_\delta\subset B_1$; an annular filling
		with matching inner trace and outer trace
		$\bn$ completes the competitor.
		Values on circles denote Sobolev traces.}
	\label{fig:fill}
\end{figure}

We translate the restriction of $\bU_\e$
to $B_\delta(a_\e)$ onto the disk $B_\delta$
centered at the origin:
\begin{equation}
	\bw_{\e,\delta}\in W^{1,p}(B_\delta;\R^2),
	\qquad
	\bw_{\e,\delta}(x)
	=\bU_\e(a_\e+x),
	\qquad x\in B_\delta.
\end{equation}
Since $\delta<1$, we have
$\overline{B_\delta}\subset B_1$.
The map $\bw_{\e,\delta}$ will define
the competitor on $B_\delta$; to complete it on
$B_1$, we seek an annular map
$\bH\in
W^{1,p}(B_1\setminus\overline{B_\delta};\R^2)$
satisfying, for every $\omega\in\Sone$,
\begin{equation}
 \begin{aligned}
 \bH_{|\partial B_\delta}(\delta\omega)
 &=(\bw_{\e,\delta})_{|\partial B_\delta}
 (\delta\omega)=\bh_{\e,\delta}(\omega),\\
 \bH_{|\partial B_1}(\omega)
 &=\bn(\omega)=\omega.
 \end{aligned}
\end{equation}
The first condition ensures that the traces match across $\partial B_\delta$, while the second enforces the prescribed boundary datum on $\partial B_1$ (where all boundary values are understood in the sense of continuous representatives of Sobolev traces). For any such $\bH$, standard Sobolev gluing yields an admissible competitor $\bV\in\cA$ defined by
\begin{equation}
	\bV(x)=
	\begin{cases}
		\bw_{\e,\delta}(x),&x\in B_\delta,\\
		\bH(x),&x\in B_1\setminus B_\delta,
	\end{cases}
\end{equation}
so that minimality of $\bu_\e$ implies $\E_\e(\bu_\e)\leqslant\E_\e(\bV)$. Notice that in the special case $\bh_{\e,\delta}(\omega)=\omega$, one could simply take $\bH=\bn$, since $\bn(\delta\omega)=\bn(\omega)=\omega$. For a general degree-one trace $\bh_{\e,\delta}$, however, a nontrivial interpolation is required; see Figure~\ref{fig:fill} for an illustration of the resulting construction.

We estimate the energy of $\bH$
using Lemma~\ref{lem:fill}, which constructs
an interpolation on an annulus with inner radius
one. To apply it to
$B_1\setminus\overline{B_\delta}$, we set
\begin{equation}
	T=\delta^{-1}>1,
	\qquad \lambda=\frac{\delta^p}{\e^2}
	=\left(\frac{\delta}{\lc}\right)^p>0,
\end{equation}
and rescale any candidate $\bH\in W^{1,p}(B_1\setminus\overline{B_\delta};\R^2)$ to a map $\bv\in W^{1,p}(B_T\setminus\overline{B_1};\R^2)$ via
\begin{equation}
	\bv(y)=\bH(\delta y),
	\qquad y\in B_T\setminus\overline{B_1}.
\end{equation}
Under this rescaling, the boundary conditions become $\bv_{|\partial B_1}(\omega)=\bh_{\e,\delta}(\omega)$ and $\bv_{|\partial B_T}(T\omega)=\omega$ for all $\omega\in\Sone$. For $\lambda>0$, an open set $D\subset\R^2$, and $\bv\in W^{1,p}_{\mathrm{loc}}(D;\R^2)$, we introduce the rescaled energy functional
\begin{equation}\label{eq:f_lambda}
	\cF_\lambda(\bv;D)
	=\int_D\left(
	\frac{|\nabla\bv|^p}{p}
	+\frac{\lambda}{4}(1-|\bv|^2)^2
	\right)\dd\in[0,\infty],
\end{equation}
so that $\cF_1=\cF$ from \eqref{eq:scale}. A change of variables yields
\begin{equation}\label{eq:scale_loc}
	\E_\e(\bH;B_1\setminus\overline{B_\delta})
	=\delta^{2-p}
	\cF_\lambda(\bv;B_T\setminus\overline{B_1}),
\end{equation}
which reduces the problem to estimating $\cF_\lambda$ on the expanding annulus $B_T\setminus\overline{B_1}$. Note that when $\delta/\lc$ is large, so is the coupling parameter $\lambda$, heavily penalizing deviations of $|\bv|$ from $1$.

\begin{figure}[t]
	\centering
	\resizebox{\linewidth}{!}{%
		\includegraphics[]{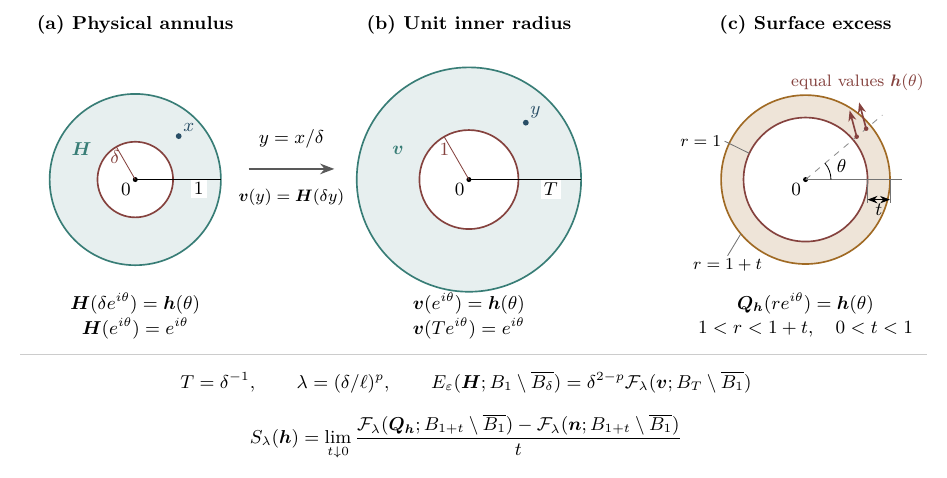}%
	}
	\caption{Rescaling and surface excess. The change
		$y=x/\delta$ gives $T=\delta^{-1}$ and
		$\lambda=(\delta/\lc)^p$.
		Panel (c) represents the radially constant map
		$\bQ_{\bh}$ on a thin
		annulus. Its parallel arrows indicate equal
		target values on the same ray. The difference
		quotient measures excess energy per unit radial
		thickness relative to $\bn$.
		The panels use separate drawing scales.}
	\label{fig:scale}
\end{figure}

We next introduce the boundary quantity that governs this interpolation energy. Setting $\TT=\R/(2\pi\mathbb Z)$ and identifying $\R^2$ with $\C$, we may write any nonvanishing degree-one trace $\bh\in W^{1,p}(\TT;\C)$ in polar form as
\begin{equation}
	\bh(\theta)
	=\rh(\theta)e^{i(\theta+\ps(\theta))},
	\qquad \theta\in\TT,
\end{equation}
for some $\rh,\ps\in W^{1,p}(\TT;\R)$ with $\rh>0$ (in our setting, $\bh(\theta)=\bh_{\e,\delta}(e^{i\theta})$). Extending $\bh$ by being a constant along each ray,
we define $\bQ_{\bh}\in W^{1,p}_{\mathrm{loc}}(\R^2\setminus\overline{B_1};\C)$
by setting
\begin{equation}
 \bQ_{\bh}(re^{i\theta})
 =\bh(\theta),
 \qquad (r,\theta)\in(1,\infty)\times\TT.
\end{equation}
Thus $\bQ_{\bh}$ retains
the same angular profile on every concentric circle. See Figure~\ref{fig:scale} for a schematic illustration. 

Next, we define the \emph{surface excess} by
\begin{equation}\label{eq:surf}
 \Surf_\lambda(\bh)
 :=\frac1p\int_0^{2\pi}|\bh'|^p\,d\theta
 +\frac{\lambda}{4}\int_0^{2\pi}(1-\rh^2)^2\,d\theta
 -\frac{2\pi}{p}.
\end{equation}
This is the boundary energy of $\bh$
relative to the radial trace, with arc-length
measure and tangential derivatives. Indeed, since $|\nabla\bQ_{\bh}(re^{i\theta})|=r^{-1}|\bh'(\theta)|$ almost everywhere, integrating in polar coordinates
\begin{equation}
 \Surf_\lambda(\bh)
 =\lim_{t\downarrow0}
 \frac{
 \cF_\lambda(\bQ_{\bh};
 B_{1+t}\setminus\overline{B_1})
 -\cF_\lambda(\bn;
 B_{1+t}\setminus\overline{B_1})
 }{t}.
\end{equation}
When $\rh\equiv1$, Jensen's inequality and the periodicity of $\ps$ yield
\begin{equation}
	\Surf_\lambda(\bh)\geqslant0,
	\qquad\text{with}\qquad
	\Surf_\lambda(\bh)=0
	\quad\Longleftrightarrow\quad
	\bh(\theta)=e^{i\am}e^{i\theta}
	\text{ for some }\am\in\R.
\end{equation}
For $\rh\not\equiv1$, however, $\Surf_\lambda(\bh)$ may be negative, as the reduction in angular gradient energy can outweigh the potential penalty.

The filling lemma below (Lemma~\ref{lem:fill}) bounds the interpolation
energy in terms of $\Surf_\lambda(\bh)$,
under the normalization
$\langle\ps\rangle=0$.
For a general inner trace, set
\begin{equation}
 \am
 :=\frac1{2\pi}\int_0^{2\pi}\ps(\theta)\,d\theta,
 \qquad
 \widehat{\ps}:=\ps-\am,
 \qquad
 \widehat{\bh}:=e^{-i\am}\bh
 =\rh e^{i(\theta+\widehat{\ps})}.
\end{equation}
Then $\langle\widehat{\ps}\rangle=0$,
and the constant target rotation preserves the
modulus and the norm of the tangential derivative.
Consequently,
\begin{equation}
 \Surf_\lambda(\widehat{\bh})
 \overset{\eqref{eq:surf}}{=}
 \Surf_\lambda(\bh).
\end{equation}
The filling constructed for $\widehat{\bh}$
has inner boundary values
$e^{-i\am}\bh$, whereas
$\bw_{\e,\delta}$ has boundary values
$\bh$. To match these traces, we replace
the map on $B_\delta$ by
$e^{-i\am}\bw_{\e,\delta}$.
This constant rotation preserves its energy.
The two pieces can then be glued, and the filling's
outer trace remains $\bn$, so the resulting
competitor belongs to $\cA$.

\begin{lemma}[Strict exterior filling]\label{lem:fill}
	Let $p>2$, let $\TT=\R/(2\pi\mathbb Z)$,
	and let $\rh,\ps\in
	W^{1,p}(\TT;\R)$ satisfy
	$0<\rh\leqslant1$ and
	$\langle\ps\rangle=0$.
	Define the prescribed inner trace by
	\begin{equation}
		\bh\in W^{1,p}(\TT;\C),
		\qquad
		\bh(\theta)
		=\rh(\theta)e^{i(\theta+\ps(\theta))},
		\qquad \theta\in\TT.
	\end{equation}
	There exist constants
	\begin{equation}
		\eta\in(0,1/2),\qquad \lambda_0>1,\qquad T_0>2,
		\qquad C>0,\qquad
		\frac1{p-1}<k<\frac1{p-2},
	\end{equation}
	depending only on $p$, such that, whenever
	$\lambda\geqslant\lambda_0$, $T\geqslant T_0$,
	and $\rh\geqslant1-\eta$, there exists
	$\bv\in W^{1,p}(A_{1,T};\C)$,
	where $A_{1,T}=B_T\setminus\overline{B_1}$,
	with the following properties:
	\begin{enumerate}
		\item \emph{Boundary values.} For every
		$\theta\in\TT$,
		\begin{equation}
			\bv_{|\partial B_1}(e^{i\theta})=\bh(\theta),
			\qquad
			\bv_{|\partial B_T}(Te^{i\theta})=e^{i\theta}.
		\end{equation}
		
		\item \emph{Modulus bounds.} For every
		$(r,\theta)\in(1,T)\times\TT$,
		\begin{equation}
			\rh(\theta)\leqslant
			|\bv(re^{i\theta})|\leqslant1.
		\end{equation}
		
		\item \emph{Energy estimate.}
		\begin{equation}\label{eq:fill}
			\cF_\lambda(\bv;A_{1,T})
			\leqslant\cp(1-T^{2-p})
			+k\Surf_\lambda(\bh)+C\lambda^{-1/p}.
		\end{equation}
	\end{enumerate}
	If $\rh\equiv1$, the map can be chosen with
	$|\bv|\equiv1$ and the final error term
	in \eqref{eq:fill} can be omitted.
\end{lemma}

\begin{remark}[Geometry and role of the strict filling]
	The boundary conditions permit gluing to an
	interior map with trace $\bh$ and extending
	by $\bn$ outside $B_T$.
	The filling joins the inner values
	$\rh(\theta)e^{i(\theta+\ps(\theta))}$
	to $e^{i\theta}$ at the outer boundary, while
	\begin{equation}
		0\leqslant1-|\bv(re^{i\theta})|
		\leqslant1-\rh(\theta)\leqslant\eta,
		\qquad (r,\theta)\in(1,T)\times\TT.
	\end{equation}
	Thus it introduces no larger modulus defect along
	any ray. Its continuous normalized map
	$\bv/|\bv|$ has degree one on every
	intermediate circle.
	
	Since the radial energy is
	$\cF_\lambda(\bn;A_{1,T})
	=\cp(1-T^{2-p})$,
	\eqref{eq:fill} bounds the additional energy
	needed to accommodate $\bh$.
	For $\rh\equiv1$, the potential energy and
	the error term both vanish. Comparing with the
	radially constant extension then gives
	\begin{equation}
		\begin{aligned}
			\cF_\lambda(\bQ_{\bh};A_{1,T})
			-\cF_\lambda(\bn;A_{1,T})
			&=\frac{1-T^{2-p}}{p-2}\Surf_\lambda(\bh),\\
			\cF_\lambda(\bQ_{\bh};A_{1,T})
			-\cF_\lambda(\bv;A_{1,T})
			&\overset{\eqref{eq:fill}}{\geqslant}
			\left(\frac{1-T^{2-p}}{p-2}-k\right)
			\Surf_\lambda(\bh).
		\end{aligned}
	\end{equation}
	For sufficiently large $T$, the coefficient
	is positive because $k<1/(p-2)$.
	If $\Surf_\lambda(\bh)>0$,
	the filling therefore lowers the energy while
	attaining the required outer trace: the angular
	energy saved exceeds the cost of radial variation.
	For variable modulus, the estimate includes
	the error $C\lambda^{-1/p}\to0$ as
	$\lambda\to\infty$.
	
	In Proposition~\ref{prop:excess}, the same strict
	inequality yields
	$\qexp=1/k\in(p-2,p-1)$ and, for every
	fixed admissible radius $\delta>0$, the factor
	\begin{equation}
		\lc^{\qexp-(p-2)}\delta^{-\qexp}
		=\e^{2(\qexp-p+2)/p}\delta^{-\qexp}
		\longrightarrow0
		\qquad\text{as }\e\downarrow0.
	\end{equation}
	This is why the strict gain in the filling estimate
	makes the exterior energy excess asymptotically
	nonpositive at fixed radii.
\end{remark}

\begin{proof}[Proof of Lemma~\ref{lem:fill}]
We first let the phase correction decay with the radius.
The gain in tangential energy is linear in the decay
exponent, whereas the radial variation costs quadratically.
For a nonunit trace, we first restore the modulus in
a thin shell, using part of the surface potential energy.

Fix $p>2$ and
$\ps\in W^{1,p}(\TT;\R)$ with
$\langle\ps\rangle=0$.
For $T>1$, write
$A_{1,T}:=B_T\setminus\overline{B_1}$.
Given $\ad\in(0,1]$, define
$\qcut_{\ad,T}\in C^1([1,T];[0,1])$ and
$\bv\in W^{1,p}(A_{1,T};\Sone)$ by
\begin{equation}\label{eq:phase_fill}
 \qcut_{\ad,T}(r)
 :=\frac{r^{-\ad}-T^{-\ad}}
          {1-T^{-\ad}},
 \qquad
 \bv(re^{i\theta})
 :=e^{i(\theta+\qcut_{\ad,T}(r)
                    \ps(\theta))},
 \qquad (r,\theta)\in(1,T)\times\TT.
\end{equation}
Set
\begin{equation}\label{eq:phase_excess}
 e:=\Dp_p(\ps)
 :=\frac1p\int_0^{2\pi}
       \bigl(|1+\ps'|^p-1\bigr)\,d\theta.
\end{equation}
The scalar remainder $R:\R\to[0,\infty)$,
$R(s):=|1+s|^p-1-ps$, satisfies
\begin{equation}\label{eq:rem}
 c_p^\ast(s^2+|s|^p)
 \leqslant R(s)\leqslant C_p(s^2+|s|^p),
 \qquad s\in\R,
\end{equation}
for constants $c_p^\ast,C_p\in(0,\infty)$ depending only
on $p$. Indeed, strict convexity gives $R(s)>0$ for
$s\ne0$, while Taylor expansion and the behavior at
infinity give
\begin{equation}\label{eq:rem_lim}
 \lim_{s\to0}\frac{R(s)}{s^2}=\frac{p(p-1)}2,
 \qquad
 \lim_{|s|\to\infty}\frac{R(s)}{|s|^p}=1.
\end{equation}
Continuity handles the remaining compact sets.
Since $\int_0^{2\pi}\ps'=0$,
we have $pe=\int_0^{2\pi}R(\ps')$.
The zero-mean Poincar\'e inequalities therefore yield
\begin{equation}\label{eq:phase_ctrl}
 \|\ps\|_{W^{1,2}(\TT)}^2
 +\|\ps\|_{W^{1,p}(\TT)}^p
 \leqslant C_p e.
\end{equation}
In particular, $e\geqslant0$. We shall also use
\begin{equation}\label{eq:growth}
 |1+s|^p\leqslant C_p(1+R(s)),
 \qquad s\in\R.
\end{equation}

Suppose first that $\rh\equiv1$.
The map \eqref{eq:phase_fill} has the required traces
and unit modulus. With $b:=(1-T^{-\ad})^{-1}>1$,
its cutoff satisfies
\begin{equation}\label{eq:cutoff}
 0\leqslant\qcut_{\ad,T}\leqslant1,
 \qquad
 \qcut_{\ad,T}(r)\leqslant br^{-\ad},
 \qquad
 |r\qcut_{\ad,T}'(r)|
 =b\ad r^{-\ad}.
\end{equation}
Convexity gives
\begin{equation}\label{eq:ang_gain}
 \frac1p\int_0^{2\pi}
 \bigl(|1+\qcut_{\ad,T}\ps'|^p-1\bigr)
 \,d\theta
 \leqslant\qcut_{\ad,T}e.
\end{equation}
To estimate the radial contribution, use
\begin{equation}\label{eq:power_inc}
 (a^2+d^2)^{p/2}-|a|^p
 \leqslant C_p\bigl(|a|^{p-2}d^2+|d|^p\bigr),
 \qquad a,d\in\R.
\end{equation}
For $|d|\leqslant|a|$, this follows from the mean value
theorem; otherwise both terms on the left are bounded
by a multiple of $|d|^p$.
Apply it with
$a=1+\qcut_{\ad,T}\ps'$
and $d=r\qcut_{\ad,T}'\ps$.
H\"older's inequality and \eqref{eq:phase_ctrl} give
\begin{equation}\label{eq:mixed}
 \int_0^{2\pi}
 |\ps'|^{p-2}|\ps|^2\,d\theta
 \leqslant
 \|\ps'\|_{L^p(\TT)}^{p-2}
 \|\ps\|_{L^p(\TT)}^2
 \overset{\eqref{eq:phase_ctrl}}{\leqslant}C_p e.
\end{equation}
Using \eqref{eq:cutoff}, the radial
contribution integrated in $\theta$ is therefore at most
$C_pb^p\ad^2r^{-2\ad}e$.
Here $b,r\geqslant1$, $p>2$, and $\ad\leqslant1$.
Polar integration yields
\begin{equation}\label{eq:unit_en}
 \cF_\lambda(\bv;A_{1,T})
 -\cp(1-T^{2-p})
 \leqslant
 \left(
 \frac{b}{p-2+\ad}
 +\frac{C_pb^p\ad^2}{p-2+2\ad}
 \right)e.
\end{equation}
As $T\to\infty$, the coefficient tends to
\begin{equation}\label{eq:gain}
 K(\ad)
 :=\frac1{p-2+\ad}
   +\frac{C_p\ad^2}{p-2+2\ad}
 =\frac1{p-2}-\frac{\ad}{(p-2)^2}
   +O_p(\ad^2),
 \qquad \ad\downarrow0.
\end{equation}
Choose $\ad\in(0,1]$ small enough that
$K(\ad)<1/(p-2)$, then choose
\begin{equation}\label{eq:gain_par}
 \max\left\{K(\ad),\frac1{p-1}\right\}
 <k_0<k<\frac1{p-2}.
\end{equation}
Since $b\downarrow1$ as $T\to\infty$, there exists
$T_*>2$, depending only on $p$, such that the coefficient
in \eqref{eq:unit_en} is at most $k_0$
for every $T\geqslant T_*$.
This proves the unit-modulus estimate without an error term.

Now let $\rh\in W^{1,p}(\TT;\R)$ satisfy
$1-\eta\leqslant\rh\leqslant1$, where
$\eta\in(0,1/2)$ will be chosen below, and set
$\bh=\rh e^{i(\theta+\ps)}
\in W^{1,p}(\TT;\C)$.
For $\lambda>1$, define $A,P\in[0,\infty)$ by
\begin{equation}\label{eq:surf_parts}
 A:=\frac1p\int_0^{2\pi}|\rh'|^p\,d\theta,
 \qquad
 P:=\frac{\lambda}{4}\int_0^{2\pi}(1-\rh^2)^2\,d\theta.
\end{equation}
The identity
$|\bh'|^2=|\rh'|^2
+\rh^2|1+\ps'|^2$ and $p/2>1$ give
\begin{equation}\label{eq:surf_lb}
 \Surf_\lambda(\bh)
 \geqslant A+e+P
 -\frac1p\int_0^{2\pi}
 (1-\rh^p)|1+\ps'|^p\,d\theta.
\end{equation}
Put $a:=1-\rh\in W^{1,p}(\TT;\R)$,
so $0\leqslant a\leqslant\eta$.
By \eqref{eq:growth}, the last
integrand, including its factor $1/p$, is at most
$C_p\eta R(\ps')+C_pa$.
Moreover,
\begin{equation}\label{eq:surf_abs}
 C_pa\leqslant\frac{\lambda a^2}{8}+C_p\lambda^{-1},
 \qquad
 P\geqslant\frac{\lambda}{4}\int_0^{2\pi}a^2\,d\theta.
\end{equation}
Thus
\begin{equation}\label{eq:surf_ctrl}
 \Surf_\lambda(\bh)
 \geqslant(1-C_p\eta)e+A+\frac12P-C_p\lambda^{-1}.
\end{equation}
Retaining half of $P$ makes the constant error
$O_p(\lambda^{-1})$, rather than $O_p(\eta)$.

Set $t:=\lambda^{-1/p}\in(0,1)$ and $L:=1+t\in(1,2)$.
Define $m\in W^{1,p}((1,L)\times\TT;\R)$ and
$\bv\in W^{1,p}(B_L\setminus\overline{B_1};\C)$ by
\begin{equation}\label{eq:shell}
 m(r,\theta):=
 1-\left(1-\frac{r-1}{t}\right)a(\theta),
 \qquad
 \bv(re^{i\theta})
 :=m(r,\theta)e^{i(\theta+\ps(\theta))}.
\end{equation}
Here $\rh\leqslant m\leqslant1$,
$|\partial_rm|=a/t$, and
$|\partial_\theta m|\leqslant|\rh'|$.
Since $1\leqslant r\leqslant2$, its gradient energy satisfies
\begin{equation}\label{eq:shell_grad}
 \frac1p\int_{B_L\setminus\overline{B_1}}
 |\nabla\bv|^p\dd
 \leqslant C_p\left(
 t^{1-p}\int_0^{2\pi}a^p\,d\theta
 +t\int_0^{2\pi}
 \bigl(|\rh'|^p+|1+\ps'|^p\bigr)\,d\theta
 \right).
\end{equation}
Since $t^{-p}=\lambda$ and $a^p\leqslant a^2$,
the first term on the right is at most $4C_ptP$.
Also, $\int_0^{2\pi}|1+\ps'|^p=2\pi+pe$.
The shell potential is at most $2tP$, because
$\rh\leqslant m\leqslant1$ and
$\int_1^{1+t}r\,dr\leqslant2t$.
Consequently,
\begin{equation}\label{eq:shell_en}
 \cF_\lambda(\bv;
 B_L\setminus\overline{B_1})
 -\cp(1-L^{2-p})
 \leqslant C_pt(1+e+A+P).
\end{equation}

On $B_T\setminus\overline{B_L}$, use the unit-modulus
construction rescaled by $L$, with inner trace
$e^{i(\theta+\ps)}$ and outer trace
$e^{i\theta}$. If $T\geqslant2T_*$, then
$T/L\geqslant T_*$, so its excess above the radial
energy is at most $L^{2-p}k_0e\leqslant k_0e$.
The two traces agree on $\partial B_L$.
Sobolev gluing therefore gives
$\bv\in W^{1,p}(A_{1,T};\C)$ with the required
boundary values and modulus bounds. Adding the energies gives
\begin{equation}\label{eq:fill_en}
 \cF_\lambda(\bv;A_{1,T})
 -\cp(1-T^{2-p})
 \leqslant(k_0+C_pt)e+C_pt(A+P)+C_pt.
\end{equation}

Take a common constant $C_p$ in
\eqref{eq:surf_ctrl} and
\eqref{eq:fill_en}.
Choose $\eta\in(0,1/2)$ so small that
$k(1-C_p\eta)>k_0$, then $\lambda_0>1$ so large that
\begin{equation}\label{eq:fill_thresh}
 k_0+C_pt\leqslant k(1-C_p\eta),
 \qquad
 C_pt\leqslant\frac{k}{2},
 \qquad
 t=\lambda^{-1/p},
 \qquad \lambda\geqslant\lambda_0.
\end{equation}
Set $T_0:=2T_*$.
Combining \eqref{eq:fill_en},
\eqref{eq:fill_thresh}, and
\eqref{eq:surf_ctrl} yields
\begin{equation}\label{eq:fill_bound}
 \cF_\lambda(\bv;A_{1,T})
 -\cp(1-T^{2-p})
 \leqslant k\Surf_\lambda(\bh)
           +C_p(t+\lambda^{-1})
 \leqslant k\Surf_\lambda(\bh)
           +C_p\lambda^{-1/p}.
\end{equation}
All choices depend only on $p$.
When $\rh\equiv1$, use the first construction directly:
$e=\Surf_\lambda(\bh)\geqslant0$ and
$k_0<k$ give the estimate without the error term.
\end{proof}

\begin{remark}[Application to the minimizer's trace]
Let $\bu_\e\in\cA$ be a global minimizer, let
$a_\e\in B_1$ satisfy $\bu_\e(a_\e)=0$, and let
$\eta,\lambda_0,T_0$ be the constants from
Lemma~\ref{lem:fill}.
For $\e<\e_{p,\eta}$ sufficiently small, choose
\begin{equation}
 \max\{C_{p,\eta},\lambda_0^{1/p}\}\lc
 <\delta\leqslant T_0^{-1},
 \qquad \lc=\e^{2/p},
\end{equation}
such that
$\bh_{\e,\delta}(\omega)=\bU_\e(a_\e+\delta\omega)$
belongs to $W^{1,p}(\Sone;\C)$.
The interval is nonempty for sufficiently small $\e$,
and Sobolev slicing supplies the required regularity
for almost every radius.
Proposition~\ref{prop:conc}, with threshold $\eta$,
gives
\begin{equation}
 1-\eta\leqslant|\bh_{\e,\delta}|\leqslant1,
 \qquad
 \lambda=(\delta/\lc)^p\geqslant\lambda_0,
 \qquad T=\delta^{-1}\geqslant T_0.
\end{equation}
The normalized trace has degree one, as established
above using continuity of $\bU_\e/|\bU_\e|$ and the
prescribed exterior map $\bU_\e=\bn$.
Consequently, in angular coordinates we can write
\begin{equation}
 \begin{gathered}
 \bh(\theta):=\bh_{\e,\delta}(e^{i\theta})
 =\rh(\theta)e^{i(\theta+\ps(\theta))},
 \qquad \theta\in\TT,\\
 \bh\in W^{1,p}(\TT;\C),
 \qquad \rh,\ps\in W^{1,p}(\TT;\R).
 \end{gathered}
\end{equation}
To impose the phase normalization, set
\begin{equation}
 \am=\langle\ps\rangle\in\R,
 \qquad
 \widehat{\bh}=e^{-i\am}\bh
 =\rh e^{i(\theta+\ps-\am)}
 \in W^{1,p}(\TT;\C).
\end{equation}
Since $\langle\ps-\am\rangle=0$ and
$\Surf_\lambda(\widehat{\bh})=\Surf_\lambda(\bh)$,
Lemma~\ref{lem:fill} applies with prescribed
inner trace $\widehat{\bh}$.
After rescaling, its filling matches the rotated
interior map
$e^{-i\am}\bw_{\e,\delta}\in W^{1,p}(B_\delta;\C)$
on $\partial B_\delta$ and has outer trace $\bn$
on $\partial B_1$. The constant target rotation
preserves the interior energy, so the resulting
glued map is an admissible competitor for $\bu_\e$.
\end{remark}

\section{The exterior excess and its decay}\label{sec:excess}
We now apply minimality to the competitor from
Section~\ref{sec:fill}, whose annular part is supplied
by Lemma~\ref{lem:fill}. Its energy bound will control
$\bU_\e$ outside $B_\delta(a_\e)$.

Fix $\e>0$, a global minimizer $\bu_\e\in\cA$,
and a zero $a_\e\in B_1$. Recall its canonical extension
$\bU_\e\in W^{1,p}_{\mathrm{loc}}(\R^2;\R^2)$
from \eqref{eq:ext}.
The definition \eqref{eq:energy} applies directly to
$\bU_\e$ on any measurable subset of $\R^2$.
Since $\bU_\e=\bn$ outside $B_1$,
\begin{equation}
 \E_\e(\bU_\e;\R^2)
 =\E_\e(\bu_\e)
  +\E_\e(\bn;\R^2\setminus\overline{B_1})
 \overset{\eqref{eq:rad_tail}}{=}
 \E_\e(\bu_\e)+\cp<\infty.
\end{equation}

To compare the exterior energy of $\bU_\e$ with that
of a radial map centered at the chosen zero $a_\e$,
introduce
\begin{equation}\label{eq:rad_a}
 \bn_{a_\e}\in C^\infty(\R^2\setminus\{a_\e\};\Sone),
 \qquad
 \bn_{a_\e}(x)=\frac{x-a_\e}{|x-a_\e|}.
\end{equation}
Since $|\bn_{a_\e}|=1$ and
$|\nabla \bn_{a_\e}(x)|=|x-a_\e|^{-1}$,
polar integration gives
\begin{equation}
 \E_\e\bigl(
 \bn_{a_\e};\R^2\setminus\overline{B_\delta(a_\e)}
 \bigr)=\cp\delta^{2-p},
 \qquad \delta>0.
\end{equation}
We therefore define the \emph{exterior energy excess
centered at $a_\e$} by
\begin{equation}\label{eq:excess}
 \cS_\e(\delta)
 :=\E_\e\bigl(
 \bU_\e;\R^2\setminus\overline{B_\delta(a_\e)}
 \bigr)-\cp\delta^{2-p},
 \qquad \delta\in(0,\infty).
\end{equation}
This quantity compares the full exterior energy
of $\bU_\e$, including its potential energy, with
the radial reference energy. Its sign records
whether the former is above or below the latter;
equality alone does not identify the two maps.
Polar integration also shows that
$\cS_\e:(0,\infty)\to\R$ is locally
absolutely continuous.

It is possible to prove the following estimate
by comparing $\bu_\e$ with the translated interior
map and the filling from Lemma~\ref{lem:fill}.
The strict gain $k<1/(p-2)$ yields the decay exponent
$\qexp=1/k>p-2$.

\begin{proposition}[Decay of the exterior energy excess]
\label{prop:excess}
Let $p>2$. For $\e>0$, let $\bu_\e\in\cA$ be a global
minimizer of $\E_\e$, and choose $a_\e\in B_1$ with
$\bu_\e(a_\e)=0$. Set $\lc=\e^{2/p}$ and let
$\cS_\e:(0,\infty)\to\R$ be the exterior
energy excess centered at $a_\e$, defined in
\eqref{eq:excess}.

There exist thresholds
\begin{equation}
 M>1,\qquad r_*\in(0,1),\qquad \e_*\in(0,1),
\end{equation}
and constants
\begin{equation}
 C>0,\qquad \qexp\in(p-2,p-1),
\end{equation}
all depending only on $p$, such that, whenever
\begin{equation}
 0<\e<\e_*,
 \qquad M\lc\leqslant\delta\leqslant r_*,
\end{equation}
the exterior energy excess satisfies
\begin{equation}\label{eq:decay}
 \cS_\e(\delta)
 \leqslant C\lc^{\qexp-(p-2)}\delta^{-\qexp}.
\end{equation}
The thresholds and constants are uniform over all
global minimizers and all choices of their zeros.
\end{proposition}

\begin{remark}[Meaning and role of the decay estimate]
Estimate~\eqref{eq:decay} can be written as
\begin{equation}
 \E_\e\bigl(
 \bU_\e;\R^2\setminus\overline{B_\delta(a_\e)}
 \bigr)
 \leqslant
 \cp\delta^{2-p}
 \left[
 1+\frac{C}{\cp}
 \left(\frac{\lc}{\delta}\right)^{\qexp-p+2}
 \right],
 \qquad M\lc\leqslant\delta\leqslant r_*.
\end{equation}
The leading term $\cp\delta^{2-p}$ is the exterior
energy of the radial comparison map $\bn_{a_\e}$
defined in \eqref{eq:rad_a}.
Since the decay argument uses $\qexp=1/k$, and $k<1/(p-2)$ by Lemma~\ref{lem:fill}, we have
\begin{equation}
 \qexp-p+2=\frac1k-(p-2)>0,
 \qquad
 \left(\frac{\lc}{\delta}\right)^{\qexp-p+2}
 \longrightarrow0
 \quad\text{whenever }\frac{\delta}{\lc}\to\infty.
\end{equation}
Geometrically, on the exterior of a circle
$\partial B_\delta(a_\e)$ with $\delta\gg\lc$,
the energy of $\bU_\e$ is bounded by the radial
reference energy with a vanishing relative error.
At the endpoint $k=1/(p-2)$, the exponent would
vanish and the error factor would remain equal to one.

For each fixed $\delta\in(0,r_*]$, the estimate applies
for all sufficiently small $\e$, because $M\lc\to0$.
Its absolute error then tends to zero, giving
\begin{equation}
 \limsup_{\e\downarrow0}\cS_\e(\delta)
 \leqslant0.
\end{equation}
Thus the exterior energy is asymptotically no greater
than the radial reference energy. This comparison
requires no prior information on the position of
$a_\e$: the extension $\bU_\e$ makes it valid even when
$B_\delta(a_\e)$ intersects $\partial B_1$.

The estimate supplies the small upper bound needed
in Proposition~\ref{prop:defect}. That proposition
converts the energy comparison into control of
$\nabla|\bU_\e|$, the potential energy of $\bU_\e$, and
the weighted difference between the phase gradients
of $\bU_\e$ and $\bn_{a_\e}$. These separate controls
will locate the chosen zero and identify the limiting
map in Section~\ref{sec:proof}.
\end{remark}

\begin{proof}[Proof of Proposition~\ref{prop:excess}]
We compare the minimizer with a translated and rotated
copy of its core, joined to the prescribed boundary values
by Lemma~\ref{lem:fill}. Cancelling the copied core energy
gives a differential inequality for the exterior excess,
which we then integrate.

Fix $p\in(2,\infty)$, and let
$\eta\in(0,1/2)$, $\lambda_0>1$, $T_0>2$, $C>0$, and
$k\in(1/(p-1),1/(p-2))$ be the constants in
Lemma~\ref{lem:fill}. Set
\begin{equation}\label{eq:q}
 \qexp:=\frac1k\in(p-2,p-1).
\end{equation}
Choose $M,r_*,\e_*\in(0,\infty)$ such that
\begin{equation}\label{eq:excess_par}
 M>\max\{1,C_{p,\eta},\lambda_0^{1/p}\},
 \qquad
 0<r_*<T_0^{-1},
 \qquad
 0<\e_*<
 \min\left\{1,\e_{p,\eta},
                 \left(\frac{r_*}{M}\right)^{p/2}\right\}.
\end{equation}
These choices depend only on $p$.

Fix $\e\in(0,\e_*)$, a global minimizer
$\bu_\e\in\cA\subset W^{1,p}(B_1;\R^2)$ of
$\E_\e$, and a zero $a_\e\in B_1$.
Recall that $\lc=\e^{2/p}$ and that
$\bU_\e\in W^{1,p}_{\mathrm{loc}}(\R^2;\R^2)$
is the canonical extension from
\eqref{eq:ext}.
Proposition~\ref{prop:conc}, the projection bound,
and $|\bn|=1$ give
\begin{equation}\label{eq:excess_mod}
 1-\eta\leqslant|\bU_\e|\leqslant1
 \qquad\text{on }
 \R^2\setminus B_{M\lc}(a_\e).
\end{equation}
Recall that $\TT=\R/(2\pi\mathbb Z)$.
For almost every $\delta\in(M\lc,r_*)$,
Sobolev slicing gives
\begin{equation}\label{eq:trace}
 \bh_\delta\in W^{1,p}(\TT;\C),
 \qquad
 \bh_\delta(\theta)
 :=\bU_\e(a_\e+\delta e^{i\theta}),
 \qquad \theta\in\TT.
\end{equation}
Here $\e$ is fixed and is suppressed in the trace notation.
The normalized trace has degree one for every radius
in this interval, independently of the slicing regularity.
Indeed, the continuous map
$\bU_\e/|\bU_\e|$ is defined throughout
the exterior of $B_\delta(a_\e)$.
Its restrictions to concentric circles therefore have
the same degree. On a circle of radius
$R>1+|a_\e|$, its restriction is
$(a_\e+Re^{i\theta})/|a_\e+Re^{i\theta}|$.
Replacing $a_\e$ by $s a_\e$, with $s\in[0,1]$, gives
a nonvanishing homotopy from $e^{i\theta}$ to this map.
Consequently,
\begin{equation}\label{eq:trace_deg}
 \deg\left(
 \frac{\bh_\delta}{|\bh_\delta|},
 \TT\right)=1.
\end{equation}
For a radius satisfying \eqref{eq:trace},
the degree-one lifting gives
\begin{equation}\label{eq:trace_phase}
 \bh_\delta(\theta)
 =\rh(\theta)e^{i(\theta+\ps(\theta))},
 \qquad
 \rh,\ps\in
 W^{1,p}(\TT;\R),
 \qquad 1-\eta\leqslant\rh\leqslant1.
\end{equation}
To justify this representation, set
$\bm g=\bh_\delta/|\bh_\delta|
\in W^{1,p}(\TT;\Sone)$.
The lower modulus bound ensures this regularity.
The identity
$\bm g'=i(\bm g\wedge\bm g')\bm g$ almost everywhere
gives an absolutely continuous phase by integrating
$\bm g\wedge\bm g'$ and choosing its initial value.
By \eqref{eq:trace_deg}, this phase
has increment $2\pi$ over $[0,2\pi]$.
Subtracting $\theta$ therefore gives the periodic
function $\ps$ in
\eqref{eq:trace_phase}.

Normalize this phase correction by setting
\begin{equation}\label{eq:phase_mean}
 \am:=\langle\ps\rangle
 =\frac1{2\pi}\int_0^{2\pi}
       \ps(\theta)\,d\theta\in\R.
\end{equation}
For the parameters
\begin{equation}\label{eq:fill_par}
 \lambda:=\left(\frac{\delta}{\lc}\right)^p
 \overset{\eqref{eq:excess_par}}{\geqslant}
 \lambda_0,
 \qquad
 T:=\delta^{-1}
 \overset{\eqref{eq:excess_par}}{\geqslant}
 T_0,
\end{equation}
Lemma~\ref{lem:fill} applies to
$e^{-i\am}\bh_\delta$.
Denote the resulting filling by
$\bv\in W^{1,p}(B_T\setminus\overline{B_1};\C)$,
and define
\begin{equation}\label{eq:competitor}
 \bV_{\e,\delta}(x):=
 \begin{cases}
 e^{-i\am}\bU_\e(a_\e+x),
     &x\in B_\delta,\\
 \bv(x/\delta),
     &x\in B_1\setminus\overline{B_\delta}.
 \end{cases}
\end{equation}
Since $\delta<r_*<1$, the copied map is defined on
$B_\delta\Subset B_1$ and belongs to
$W^{1,p}(B_\delta;\C)$.
The rescaled filling belongs to
$W^{1,p}(B_1\setminus\overline{B_\delta};\C)$.
Their traces agree on $\partial B_\delta$, and the
outer trace is $\bn$. Sobolev gluing therefore gives
$\bV_{\e,\delta}\in\cA$.
The use of $\bU_\e$ makes the copied map
well defined even when $B_\delta(a_\e)$ intersects
$\partial B_1$.

Set $A_\delta:=B_1\setminus\overline{B_\delta}$.
Translation and constant target rotation preserve both
energy terms. Minimality and the energy of the prescribed
radial extension yield
\begin{equation}\label{eq:core_cancel}
 \E_\e(\bU_\e;\R^2)
 =\E_\e(\bu_\e)+\cp
 \leqslant
 \E_\e(\bV_{\e,\delta})+\cp
 =\E_\e
       (\bU_\e;B_\delta(a_\e))
   +\E_\e(\bV_{\e,\delta};A_\delta)
   +\cp.
\end{equation}
The scaling identity \eqref{eq:scale_loc}, the filling
estimate \eqref{eq:fill}, and invariance of the surface
excess under constant target rotations give
\begin{equation}\label{eq:ann_bound}
 \E_\e(\bV_{\e,\delta};A_\delta)
 \leqslant
 \cp(\delta^{2-p}-1)
 +k\delta^{2-p}\Surf_\lambda(\bh_\delta)
 +C\lc\delta^{1-p}.
\end{equation}
Here the last term follows from
$\delta^{2-p}\lambda^{-1/p}
=\lc\delta^{1-p}$.
Cancelling the copied core energy in
\eqref{eq:core_cancel} and using
\eqref{eq:ann_bound} and
\eqref{eq:excess}, we obtain
\begin{equation}\label{eq:surf_cmp}
 \cS_\e(\delta)
 \leqslant
 k\delta^{2-p}\Surf_\lambda(\bh_\delta)
 +C\lc\delta^{1-p}.
\end{equation}

We next relate the surface excess to the derivative
of $\cS_\e\in
W^{1,1}_{\mathrm{loc}}((0,\infty);\R)$.
Polar integration in \eqref{eq:excess}, together with
$(p-2)\cp=2\pi/p$, gives
\begin{equation}\label{eq:excess_deriv}
 -\delta\cS_\e'(\delta)
 =\delta^2\int_0^{2\pi}
 \left(
 \frac{|\nabla\bU_\e(a_\e+\delta e^{i\theta})|^p}{p}
 +\frac{(1-|\bh_\delta(\theta)|^2)^2}{4\e^2}
 \right)d\theta
 -\frac{2\pi}{p}\delta^{2-p}
\end{equation}
for almost every $\delta\in(M\lc,r_*)$.
For these radii, the tangential chain rule gives
\begin{equation}\label{eq:tan_bound}
 |\bh_\delta'(\theta)|
 \leqslant
 \delta|\nabla\bU_\e(a_\e+\delta e^{i\theta})|
 \qquad\text{for a.e. }\theta\in(0,2\pi).
\end{equation}
Since $\delta^{2-p}\lambda=\delta^2/\e^2$,
\eqref{eq:surf},
\eqref{eq:excess_deriv}, and
\eqref{eq:tan_bound} imply
\begin{equation}\label{eq:surf_deriv}
 \delta^{2-p}\Surf_\lambda(\bh_\delta)
 \leqslant-\delta\cS_\e'(\delta).
\end{equation}
Combining \eqref{eq:surf_cmp}
with \eqref{eq:surf_deriv}, we conclude that
\begin{equation}\label{eq:ode}
 \cS_\e(\delta)
 \leqslant-k\delta\cS_\e'(\delta)
          +C\lc\delta^{1-p}
 \qquad\text{for a.e. }\delta\in(M\lc,r_*).
\end{equation}
Neither the surface excess nor the exterior excess
is assumed nonnegative.

To bound the initial value, set
$R:=\lc^{-1}\in(1,\infty)$.
Scaling and \eqref{eq:mono} give
\begin{equation}\label{eq:en_bound}
 \E_\e(\bU_\e;\R^2)
 =\lc^{2-p}\J(R)
 \overset{\eqref{eq:mono}}{\leqslant}
 \J(1)\lc^{2-p}.
\end{equation}
The number $\J(1)\in(0,\infty)$
depends only on $p$; its finiteness follows by using
$x\mapsto x\in W^{1,p}(B_1;\R^2)$ as a competitor.
Nonnegativity of the energy density and
\eqref{eq:en_bound} therefore imply
\begin{equation}\label{eq:excess_init}
 \cS_\e(M\lc)
 \leqslant
 \E_\e(\bU_\e;\R^2)
 \overset{\eqref{eq:en_bound}}{\leqslant}
 \J(1)\lc^{2-p}.
\end{equation}

Multiplying \eqref{eq:ode} by
$\qexp\delta^{\qexp-1}$ and using
$\qexp k=1$, we obtain, after increasing
$C\in(0,\infty)$ by a factor depending only on $p$,
\begin{equation}\label{eq:ode_weight}
 \bigl(\delta^{\qexp}
       \cS_\e(\delta)\bigr)'
 \overset{\eqref{eq:ode}}{\leqslant}
 C\lc\delta^{\qexp-p}
 \qquad\text{for a.e. }\delta\in(M\lc,r_*).
\end{equation}
Local absolute continuity permits integration from
$M\lc$ to every $\delta\in[M\lc,r_*]$.
Since $\qexp<p-1$, we have
\begin{equation}\label{eq:power_int}
 \int_{M\lc}^{\delta}
 t^{\qexp-p}\,dt
 =
 \frac{(M\lc)^{\qexp-p+1}
       -\delta^{\qexp-p+1}}
      {p-1-\qexp}
 \leqslant
 \frac{(M\lc)^{\qexp-p+1}}
      {p-1-\qexp}.
\end{equation}
Integrating \eqref{eq:ode_weight}
and applying \eqref{eq:excess_init} and
\eqref{eq:power_int} gives
\begin{equation}\label{eq:excess_bound}
 \delta^{\qexp}\cS_\e(\delta)
 \leqslant
 \left(
 \J(1)M^{\qexp}
 +\frac{CM^{\qexp-p+1}}{p-1-\qexp}
 \right)\lc^{\qexp-p+2}.
\end{equation}
Dividing by $\delta^{\qexp}>0$ proves
\eqref{eq:decay}. Every constant and threshold depends
only on $p$, so the estimate is uniform over all global
minimizers and all choices of the zero $a_\e$.
\end{proof}

\section{Convexity and the exterior defect}\label{sec:defect}
Proposition~\ref{prop:excess} compares the energy of $\bU_\e$
on $\R^2\setminus\overline{B_\delta(a_\e)}$ with
$\cp\delta^{2-p}$, the energy of $\bn_{a_\e}$
on that exterior domain.
We now turn this comparison into estimates for the
modulus $\rh=|\bU_\e|$ and the phase gradient $\bxi$,
given locally, wherever $\bU_\e\ne0$, by
$\bU_\e=\rh e^{i\vp}$ and $\bxi=\nabla\vp$.
In Section~\ref{sec:proof}, these estimates will
first show that $a_\e\to0$ and then identify the limiting
map on each fixed annulus.

For $z,w\in\R^2$, define the convexity remainder
\begin{equation}\label{eq:breg}
 \Bp(z,w)
 =\frac{|z|^p-|w|^p}{p}
  -|w|^{p-2}w\cdot(z-w).
\end{equation}
For $f:\R^2\to\R$, $f(z)=|z|^p/p$, this is the difference
between $f(z)$ and the affine approximation
$f(w)+\nabla f(w)\cdot(z-w)$, centered at $w$ and evaluated at $z$. Uniform convexity gives a constant
$b_p>0$, depending only on $p$, such that
\begin{equation}\label{eq:convex}
 \Bp(z,w)\geqslant
 b_p\bigl(|z-w|^p+|w|^{p-2}|z-w|^2\bigr)
 \qquad\text{for all }z,w\in\R^2.
\end{equation}
In particular, $\Bp$ is nonnegative and vanishes precisely
when $z=w$. We will use the two terms on the right-hand
side of \eqref{eq:convex} separately. The term
$b_p|w|^{p-2}|z-w|^2$, quadratic in $z$, will show that the chosen zero
$a_\e$ tends to the origin, whereas $b_p|z-w|^p$ will
give strong $L^p$ convergence of the weighted phase
difference defined in \eqref{eq:phase_err}.

To apply \eqref{eq:convex}, we introduce the phase gradients
of the radial comparison map and the minimizer.
For each $a\in\R^2$, define
the \emph{radial vortex phase gradient centered at $a$} by
\begin{equation}
 \bxi_a\in C^\infty(\R^2\setminus\{a\};\R^2),
 \qquad
 \bxi_a(x)=\frac{(x-a)^\perp}{|x-a|^2},
\end{equation}
where $y^\perp=(-y_2,y_1)$ for $y=(y_1,y_2)\in\R^2$.
This terminology reflects the identity
$\bxi_a=\nabla\theta_a$ wherever the radial map
$x\mapsto(x-a)/|x-a|$ is written locally as $e^{i\theta_a}$.

We next define the corresponding field for the minimizer.
Fix $\e>0$, a global minimizer $\bu_\e\in\cA$, and a zero
$a_\e\in B_1$. Let $\bU_\e$ be its extension by $\bn$ from
\eqref{eq:ext}. For $\delta>0$, consider
the exterior domain
$\Dext_{\e,\delta}=\R^2\setminus\overline{B_\delta(a_\e)}$
and assume that $|\bU_\e|\geqslant1/2$ there.
We define the modulus and phase gradient by\footnote{For
$v,w\in\R^2$, we write
$v\wedge w=v_1w_2-v_2w_1=\det(v,w)$.
For a Sobolev map $\bU=(\Ucomp_1,\Ucomp_2)$, this notation is applied
componentwise:
$\bU\wedge\nabla \bU=(\bU\wedge\partial_1\bU,\bU\wedge\partial_2\bU)
=\Ucomp_1\nabla \Ucomp_2-\Ucomp_2\nabla \Ucomp_1$.}
\begin{equation}\label{eq:phase}
 \rh=|\bU_\e|
 \in W^{1,p}_{\mathrm{loc}}(\Dext_{\e,\delta};\R),
 \qquad
 \bxi=\frac{\bU_\e\wedge\nabla \bU_\e}{\rh^2}
 \in L^p(\Dext_{\e,\delta};\R^2).
\end{equation}
To interpret this formula, on every open ball
$B\Subset \Dext_{\e,\delta}$ a Sobolev lifting gives
$\vp\in W^{1,p}(B;\R)$ such that
$\bU_\e=\rh e^{i\vp}$. Then
$\bU_\e\wedge\nabla \bU_\e=\rh^2\nabla\vp$, so
$\bxi=\nabla\vp$ on $B$. Thus $\bxi$ is a globally
defined phase gradient, although the degree-one map
has no single-valued global phase. In particular, this 
definition agrees with the radial comparison field:
for $\bn_a(x)=(x-a)/|x-a|$, direct differentiation gives
$\bn_a\wedge\nabla \bn_a=\bxi_a$ on $\R^2\setminus\{a\}$.

We compare $\bxi$, the phase gradient of $\bU_\e$, with
$\bxi_{a_\e}$, the phase gradient of $\bn_{a_\e}$.
Retaining the convexity remainder and using the degree-one
circulation of $\bU_\e$, it is possible to prove the following
estimate. It controls $\nabla\rh$, the weighted difference
between these phase gradients, and the potential energy
of $\bU_\e$ on the exterior domain, where $\rh=|\bU_\e|$.

\begin{proposition}[Exterior defect estimate]\label{prop:defect}
There exists $C_p>0$, depending only on $p$, with the
following property. Let $\e>0$, let $\bu_\e\in\cA$ be
a global minimizer of $\E_\e$, and let $a_\e\in B_1$
satisfy $\bu_\e(a_\e)=0$. Let $\bU_\e$ be the extension
defined in \eqref{eq:ext}.
For every $\delta>0$ such that
$|\bU_\e|\geqslant1/2$ on
$\Dext_{\e,\delta}=\R^2\setminus\overline{B_\delta(a_\e)}$,
one has
\begin{equation}\label{eq:defect}
 \int_{\Dext_{\e,\delta}}
 \left\{
 \frac{|\nabla\rh|^p}{p}
 +\Bp\bigl(\rh\bxi,\rh^{-1/(p-1)}\bxi_{a_\e}\bigr)
 +\frac{(1-\rh^2)^2}{8\e^2}
 \right\}\dd
 \leqslant \cS_\e(\delta)+C_p\e^2\delta^{2-2p}.
\end{equation}
Here $\rh=|\bU_\e|$ and
$\bxi=\rh^{-2}(\bU_\e\wedge\nabla \bU_\e)$ are defined
in \eqref{eq:phase}, and $\cS_\e(\delta)$ is
the exterior energy excess from \eqref{eq:excess}.
\end{proposition}

\begin{remark}[Choice of the weighted fields]
The weights are chosen to connect the phase energy
of $\bU_\e$ with its degree. Indeed, the polar decomposition
gives $|\nabla \bU_\e|^2=|\nabla\rh|^2+|\rh\bxi|^2$.
Thus $z=\rh\bxi$ is the field whose $p$-power appears
in the phase contribution to the energy.

We choose the comparison field $w$ so that the linear
term in the convexity expansion,
$|w|^{p-2}w\cdot z$, becomes independent of $\rh$.
This requires $|w|^{p-2}w
 =\rh^{-1}|\bxi_{a_\e}|^{p-2}\bxi_{a_\e}$, hence $ w=\rh^{-1/(p-1)}\bxi_{a_\e}$.
With this choice, $ |w|^{p-2}w\cdot z
 =|\bxi_{a_\e}|^{p-2}\bxi_{a_\e}\cdot\bxi$.
The integral of this expression is determined by the
degree-one circulation of $\bU_\e$ on circles centered
at $a_\e$. The remaining modulus-dependent error can
be absorbed into half of the potential energy, as
explained in the proof below. Thus the weights are
dictated by the phase energy and the circulation identity.
\end{remark}

\begin{remark}
The usefulness of \eqref{eq:defect} is that the vanishing
upper bound for the exterior excess forces three separate errors to vanish:
variation of the modulus $\rh=|\bU_\e|$, the difference
between the weighted phase fields $\rh\bxi$ and
$\rh^{-1/(p-1)}\bxi_{a_\e}$, and the potential energy
of $\bU_\e$. Indeed, fix $\delta\in(0,r_*]$.
Proposition~\ref{prop:conc}, with $\eta=1/2$,
ensures $|\bU_\e|\geqslant1/2$ on $\Dext_{\e,\delta}$
for sufficiently small $\e$, so \eqref{eq:defect}
applies. Proposition~\ref{prop:excess} bounds its
right-hand side by
\begin{equation}
 C\lc^{\qexp-(p-2)}\delta^{-\qexp}
 +C_p\e^2\delta^{2-2p}\longrightarrow0.
\end{equation}
Since the three integrals on the left of
\eqref{eq:defect} are nonnegative, each tends to zero.
This provides the separate estimates needed to locate
$a_\e$ and identify the limiting map in
Section~\ref{sec:proof}.

In Section~\ref{sec:proof}, we will use twice the fact
that, for every fixed $\delta\in(0,r_*]$,
\begin{equation}
 \int_{\Dext_{\e,\delta}}
 \Bp\bigl(\rh\bxi,\rh^{-1/(p-1)}\bxi_{a_\e}\bigr)\dd
 \longrightarrow0
 \qquad\text{as }\e\downarrow0.
\end{equation}
First, outside $B_1$, the prescribed extension satisfies
$\bU_\e=\bn$, so $\rh=1$, $\bxi=\bxi_0$, and the convexity
remainder becomes $\Bp(\bxi_0,\bxi_{a_\e})$.
Its quadratic lower bound in \eqref{eq:convex},
integrated over a fixed exterior annulus, will give
$a_\e\to0$ in Step~1. Concentration then yields
$\rh\to1$ uniformly on each fixed annulus $A_r$.
Second, in Step~3, the $p$-power lower bound in
\eqref{eq:convex} will control the weighted phase
difference on $A_r$. Together with $a_\e\to0$ and
$\rh\to1$, this will give
$\bxi\to\bxi_0$ in $L^p(A_r;\R^2)$.
The other two terms in \eqref{eq:defect} will give
$\nabla\rh\to0$ in $L^p(A_r;\R^2)$ and vanishing
potential energy of $\bU_\e$ on $A_r$.
\end{remark}

\begin{proof}[Proof of Proposition~\ref{prop:defect}]
We retain the convexity remainder in the phase energy and
evaluate its linear term through the degree-one circulation.
The resulting modulus error is absorbed into half of the
potential energy.

Fix $p\in(2,\infty)$, $\e,\delta\in(0,\infty)$, a global
minimizer $\bu_\e\in\cA\subset W^{1,p}(B_1;\R^2)$
of $\E_\e$, and a zero $a_\e\in B_1$, under the
hypotheses of the proposition.
Let $\bU_\e\in W^{1,p}_{\mathrm{loc}}(\R^2;\R^2)$
be the extension from \eqref{eq:ext}.
On $\Dext_{\e,\delta}
=\R^2\setminus\overline{B_\delta(a_\e)}$, set
$\rh=|\bU_\e|
\in W^{1,p}_{\mathrm{loc}}(\Dext_{\e,\delta};\R)$.
The projection argument in Section~\ref{sec:conc},
together with the hypothesis of the proposition, gives
\begin{equation}\label{eq:def_mod}
 \frac12\leqslant\rh\leqslant1
 \qquad\text{on }\Dext_{\e,\delta}.
\end{equation}
Finite total energy implies
$\nabla\bU_\e\in L^p(\R^2;\R^{2\times2})$.
The Sobolev chain rule and \eqref{eq:phase} give
$\nabla\rh,\bxi
\in L^p(\Dext_{\e,\delta};\R^2)$, with
\begin{equation}\label{eq:chain}
 |\nabla\rh|\leqslant|\nabla\bU_\e|,
 \qquad
 |\bxi|\leqslant
 \rh^{-1}|\nabla\bU_\e|
 \overset{\eqref{eq:def_mod}}{\leqslant}
 2|\nabla\bU_\e|
\end{equation}
almost everywhere on $\Dext_{\e,\delta}$.

To obtain the polar decomposition, define
$\bm{\nu}=\bU_\e/\rh
\in W^{1,p}_{\mathrm{loc}}(\Dext_{\e,\delta};\Sone)$.
Its gradient belongs to
$L^p(\Dext_{\e,\delta};\R^{2\times2})$, and we use
its continuous representative.
For $j\in\{1,2\}$, the identity $|\bm{\nu}|^2=1$ gives
$\bm{\nu}\cdot\partial_j\bm{\nu}=0$, hence
$\partial_j\bm{\nu}
=(\bm{\nu}\wedge\partial_j\bm{\nu})\bm{\nu}^{\perp}
=(\bxi)_j\bm{\nu}^{\perp}$ almost everywhere.
Consequently,
\begin{equation}\label{eq:polar_deriv}
 \partial_j\bU_\e
 =(\partial_j\rh)\bm{\nu}
  +\rh(\bxi)_j\bm{\nu}^{\perp},
 \qquad j\in\{1,2\}.
\end{equation}
Since $\bm{\nu}$ and $\bm{\nu}^{\perp}$ are orthonormal
and $p/2>1$, we obtain
\begin{equation}\label{eq:polar_en}
 |\nabla\bU_\e|^p
 \overset{\eqref{eq:polar_deriv}}{=}
 \bigl(|\nabla\rh|^2+\rh^2|\bxi|^2\bigr)^{p/2}
 \geqslant|\nabla\rh|^p+|\rh\bxi|^p
\end{equation}
almost everywhere on the exterior domain.

We next justify the circulation identity. Define the
positively oriented unit tangent field by
\begin{equation}\label{eq:tau}
 \btau\in C^\infty(\R^2\setminus\{a_\e\};\Sone),
 \qquad
 \btau(x)=\frac{(x-a_\e)^\perp}{|x-a_\e|},
 \qquad x\ne a_\e.
\end{equation}
Identifying $\R^2$ with $\C$, for every $t>\delta$ the
continuous loop $\theta\mapsto\bm{\nu}(a_\e+t e^{i\theta})$,
$\theta\in[0,2\pi]$, has degree one.
Indeed, choose $R>\max\{t,1+|a_\e|\}$ and increase the
radius from $t$ to $R$. The normalized map remains continuous
throughout this annulus. On the outer circle it equals
$(a_\e+Re^{i\theta})/|a_\e+Re^{i\theta}|$.
Replacing $a_\e$ by $s a_\e$, with $s\in[0,1]$, gives
a homotopy from $e^{i\theta}$ to this loop; its denominator
never vanishes because $R>|a_\e|$.

For almost every $t>\delta$, Sobolev slicing gives
$\bm g\in W^{1,p}(0,2\pi;\Sone)$, defined by
$\bm g(\theta)=\bm{\nu}(a_\e+t e^{i\theta})$, with
$\bm g(0)=\bm g(2\pi)$. Its derivative agrees with the
tangential derivative of $\bm{\nu}$ along this circle.
Choose $\beta_0\in\R$ such that $e^{i\beta_0}=\bm g(0)$
and define
\begin{equation}\label{eq:lift}
 \beta\in W^{1,p}(0,2\pi;\R),
 \qquad
 \beta(\theta)=\beta_0+
 \int_0^\theta\bm g(s)\wedge\bm g'(s)\,ds,
 \qquad \theta\in[0,2\pi].
\end{equation}
Since $\bm g'=i\bm g\,(\bm g\wedge\bm g')$ almost
everywhere, $(e^{-i\beta}\bm g)'=0$, and the choice of
$\beta_0$ yields $\bm g=e^{i\beta}$.
The degree-one property therefore gives
$\beta(2\pi)-\beta(0)=2\pi$.
Using the tangential derivative identity
\begin{equation}\label{eq:slice_deriv}
 \bm g(\theta)\wedge\bm g'(\theta)
 =t\,\bxi(a_\e+t e^{i\theta})
      \cdot\btau(a_\e+t e^{i\theta})
 \quad\text{for a.e. }\theta\in(0,2\pi),
\end{equation}
we conclude that
\begin{equation}\label{eq:circ}
 \int_{\partial B_t(a_\e)}
 \bxi\cdot\btau\,d\mathcal H^1
 \overset{\eqref{eq:slice_deriv}}{=}
 \int_0^{2\pi}\bm g\wedge\bm g'\,d\theta
 \overset{\eqref{eq:lift}}{=}
 \beta(2\pi)-\beta(0)=2\pi
\end{equation}
for almost every $t>\delta$.

The reference phase field
$\bxi_{a_\e}\in
C^\infty(\R^2\setminus\{a_\e\};\R^2)$ satisfies
$\bxi_{a_\e}(x)=\btau(x)/|x-a_\e|$.
In particular, its restriction to $\Dext_{\e,\delta}$
belongs to $L^p(\Dext_{\e,\delta};\R^2)$, since
\begin{equation}\label{eq:rad_int}
 \int_{\Dext_{\e,\delta}}
 |\bxi_{a_\e}|^p\dd
 =2\pi\int_\delta^\infty t^{1-p}\,dt
 =\frac{2\pi}{p-2}\delta^{2-p}<\infty.
\end{equation}
Thus $|\bxi_{a_\e}|^{p-1}\in
L^{p/(p-1)}(\Dext_{\e,\delta};\R)$.
H\"older's inequality makes the following integral absolutely
convergent, and polar integration yields
\begin{equation}\label{eq:circ_int}
 \int_{\Dext_{\e,\delta}}
 \bxi\cdot
 |\bxi_{a_\e}|^{p-2}\bxi_{a_\e}\dd
 \overset{\eqref{eq:circ}}{=}
 2\pi\int_\delta^\infty t^{1-p}\,dt
 =\frac{2\pi}{p-2}\delta^{2-p}.
\end{equation}

We now apply \eqref{eq:breg} pointwise to the weighted fields
\begin{equation}\label{eq:weights}
 \bm z=\rh\bxi,
 \qquad
 \bm w=\rh^{-1/(p-1)}\bxi_{a_\e},
 \qquad
 \bm z,\bm w\in L^p(\Dext_{\e,\delta};\R^2).
\end{equation}
The powers of $\rh$ cancel in
$|\bm w|^{p-2}\bm w\cdot\bm z$, giving
\begin{equation}\label{eq:phase_en}
 \frac{|\bm z|^p}{p}
 \overset{\eqref{eq:breg}}{=}
 \Bp(\bm z,\bm w)
 +\bxi\cdot|\bxi_{a_\e}|^{p-2}
                         \bxi_{a_\e}
 -\frac{p-1}{p}\rh^{-p/(p-1)}
                         |\bxi_{a_\e}|^p.
\end{equation}
Every term is integrable: the mixed term is controlled by
H\"older's inequality, and \eqref{eq:breg} gives
$\Bp(\bm z,\bm w)
\in L^1(\Dext_{\e,\delta};\R)$.
By \eqref{eq:circ_int} and
\eqref{eq:rad_int}, the linear term
and the part of the last term with $\rh$ replaced
by $1$ contribute
\begin{equation}\label{eq:rad_term}
 \left(1-\frac{p-1}{p}\right)
 \frac{2\pi}{p-2}\delta^{2-p}
 =\cp\delta^{2-p}.
\end{equation}

The remaining modulus error has a fixed sign.
The mean value theorem on $[1/2,1]$ gives a constant
$K_p\in(0,\infty)$, depending only on $p$, such that
\begin{equation}\label{eq:mod_err}
 0\leqslant\frac{p-1}{p}
 \bigl(\rh^{-p/(p-1)}-1\bigr)
 |\bxi_{a_\e}|^p
 \leqslant K_p(1-\rh^2)|\bxi_{a_\e}|^p.
\end{equation}
One may take $K_p=2^{1+p/(p-1)}$: bound the absolute
derivative of $s\mapsto s^{-p/(p-1)}$ on $[1/2,1]$ and
use $1-s\leqslant1-s^2$ there. Completing the square gives
\begin{equation}\label{eq:young}
 K_p(1-\rh^2)|\bxi_{a_\e}|^p
 \leqslant\frac{(1-\rh^2)^2}{8\e^2}
       +2K_p^2\e^2|\bxi_{a_\e}|^{2p}.
\end{equation}
These inequalities hold almost everywhere on the exterior
domain. Combining \eqref{eq:polar_en},
\eqref{eq:phase_en}, and
\eqref{eq:rad_term}, and absorbing
the error by
\eqref{eq:mod_err}--\eqref{eq:young},
we obtain
\begin{equation}\label{eq:def_bound}
 \begin{split}
 &\int_{\Dext_{\e,\delta}}
 \left\{\frac{|\nabla\rh|^p}{p}
 +\Bp(\bm z,\bm w)
 +\frac{(1-\rh^2)^2}{8\e^2}\right\}\dd\\
 &\qquad\overset{\eqref{eq:young}}{\leqslant}
 \E_\e(\bU_\e;\Dext_{\e,\delta})
 -\cp\delta^{2-p}
 +2K_p^2\e^2\int_{\Dext_{\e,\delta}}
                   |\bxi_{a_\e}|^{2p}\dd.
 \end{split}
\end{equation}
Finally, polar integration gives
\begin{equation}\label{eq:tail_err}
 \int_{\Dext_{\e,\delta}}
 |\bxi_{a_\e}|^{2p}\dd
 =2\pi\int_\delta^\infty t^{1-2p}\,dt
 =\frac{\pi}{p-1}\delta^{2-2p}.
\end{equation}
Substituting \eqref{eq:tail_err} into
\eqref{eq:def_bound} and using
\eqref{eq:excess} and \eqref{eq:weights}
proves \eqref{eq:defect} with
$C_p=2\pi K_p^2/(p-1)\in(0,\infty)$, depending only on $p$.
\end{proof}

\section{Proof of the main theorem}\label{sec:proof}
We combine Propositions~\ref{prop:conc},
\ref{prop:excess}, and \ref{prop:defect}.
Proposition~\ref{prop:conc} confines every fixed modulus-defect
set to a small ball around a zero, but leaves the position of that
ball undetermined. We first locate it by combining the exterior
excess estimate of Proposition~\ref{prop:excess} with the defect
estimate of Proposition~\ref{prop:defect}. The prescribed extension
outside $B_1$ will force its center to approach the origin.
The same estimates will then control the modulus and phase on each
fixed annulus. Finally, the boundary trace will determine the map
itself from these gradient estimates.

\begin{proof}[Proof of Theorem~\ref{thm:main}]

Fix $p\in(2,\infty)$ and an arbitrary family of global minimizers
$(\bu_\e)_{\e>0}\subset\cA$, where
$\cA\subset W^{1,p}(B_1;\R^2)$ is the class in \eqref{eq:class}.
For each $\e>0$, choose any zero $a_\e\in B_1$ of $\bu_\e$;
the existence of such a zero was proved in Section~\ref{sec:conc}.
Recall the natural scale $\lc=\e^{2/p}$ and the extension
$\bU_\e\in W^{1,p}_{\mathrm{loc}}(\R^2;\R^2)$ from \eqref{eq:ext}:
\begin{equation}\label{eq:proof_ext}
 \bU_\e(x)=
 \begin{cases}
  \bu_\e(x),&x\in B_1,\\
  \bn(x),&x\in\R^2\setminus B_1,
 \end{cases}
 \qquad \bn(x)=\frac{x}{|x|}\quad(x\ne0).
\end{equation}
We use continuous representatives throughout. The projection
argument in Section~\ref{sec:conc} and $|\bn|=1$ give
\begin{equation}\label{eq:proof_modulus}
 \rh_\e:=|\bU_\e|\in W^{1,p}_{\mathrm{loc}}(\R^2;\R),
 \qquad 0\leqslant\rh_\e\leqslant1,
 \qquad \rh_\e=1\quad\text{on }\R^2\setminus B_1.
\end{equation}
We proceed in four steps, keeping the original family of minimizers
throughout the argument.

\smallskip\noindent
\emph{Step 1: the core center tends to zero.}
We claim that the chosen zeros satisfy $a_\e\to0$ as
$\e\downarrow0$. We first record the estimate that will also be used in Step~3.
For $\delta>0$, recall the exterior domain and its energy excess:
\begin{equation}\label{eq:proof_exterior}
 \Dext_{\e,\delta}:=\R^2\setminus\overline{B_\delta(a_\e)},
 \qquad
 \cS_\e(\delta)
 \overset{\eqref{eq:excess}}{=}
 \E_\e(\bU_\e;\Dext_{\e,\delta})-\cp\delta^{2-p}.
\end{equation}
Here $\cp=2\pi/[p(p-2)]>0$, and the subtracted term is the
energy of the radial map centered at $a_\e$ on the same domain.
Let $M>1$, $r_*\in(0,1)$, $\e_*\in(0,1)$, $C>0$, and
$\qexp\in(p-2,p-1)$ be the constants in
Proposition~\ref{prop:excess}. Its estimate \eqref{eq:decay}, i.e., $\cS_\e(\delta)\leqslant C\lc^{\qexp-(p-2)}\delta^{-\qexp}$ requires
$\e<\e_*$ and $M\lc\leqslant\delta\leqslant r_*$.
To apply Proposition~\ref{prop:defect} as well, we need
$\rh_\e\geqslant1/2$ on $\Dext_{\e,\delta_0}$.

Fix any $\delta_0\in(0,r_*]$, independently of $\e$.
All these requirements hold once
\begin{equation}\label{eq:proof_range}
 0<\e<\min\{\e_*,\e_{p,1/2}\},
 \qquad \max\{M,C_{p,1/2}\}\lc<\delta_0.
\end{equation}
Indeed, Proposition~\ref{prop:conc} with threshold $\eta=1/2$
gives
\begin{equation}\label{eq:proof_defect_set}
 \{x\in B_1:|\bu_\e(x)|\leqslant1/2\}
 \overset{\eqref{eq:conc}}{\subset}
 B_{C_{p,1/2}\lc}(a_\e)
 \overset{\eqref{eq:proof_range}}{\subset}B_{\delta_0}(a_\e).
\end{equation}
Thus $|\bu_\e|>1/2$ on $B_1\cap\Dext_{\e,\delta_0}$.
On the remaining part of $\Dext_{\e,\delta_0}$, the extension has
modulus one by \eqref{eq:proof_modulus}. Hence
$1/2\leqslant\rh_\e\leqslant1$ throughout this exterior domain.
Since $\lc\to0$, condition~\eqref{eq:proof_range} is satisfied
for all sufficiently small $\e$ at every fixed such $\delta_0$.

Letting $\bxi_\e$ denote the phase gradient defined in \eqref{eq:phase} and retaining the convexity remainder from \eqref{eq:defect}, we obtain
\begin{equation}\label{eq:tempmathcalB}
 \int_{\Dext_{\e,\delta_0}}
 \Bp\bigl(\rh_\e\bxi_\e,
          \rh_\e^{-1/(p-1)}\bxi_{a_\e}\bigr)\dd
 \leqslant \cS_\e(\delta_0)+C_p\e^2\delta_0^{2-2p}
 \leqslant \Herr_\e(\delta_0),
\end{equation}
where the \emph{combined error function} $\Herr_\e:(0,\infty)\to(0,\infty)$ is given by
\begin{equation}\label{eq:err}
 \Herr_\e(\delta)
 :=C\lc^{\qexp-(p-2)}\delta^{-\qexp}
   +C_p\e^2\delta^{2-2p},
 \qquad \delta>0,
\end{equation}
with $C>0$ and $C_p>0$ denoting the constants from Propositions~\ref{prop:excess} and~\ref{prop:defect}, respectively.
To locate $a_\e$, we restrict the defect estimate to a
region where $\bU_\e$ is prescribed. Choose the fixed
annulus
\begin{equation}
 K:=\{x\in\R^2:2<|x|<3\}
   =B_3\setminus\overline{B_2}.
\end{equation}
On $K$, the identity $\bU_\e=\bn$ gives
$\rh_\e=1$ and $\bxi_\e=\bxi_0$, so the defect compares
the explicit radial phase fields centered at $0$ and
at $a_\e$.\footnote{Here $\bxi_0$ denotes the radial
phase field with center $a=0$, not $\bxi_\e$ evaluated
at $\e=0$. The radii $2$ and $3$ are chosen for
convenience: $K$ lies outside every excluded disk
$B_{\delta_0}(a_\e)$ with $a_\e\in B_1$ and
$\delta_0<1$, while its boundedness and fixed positive
area will give a lower bound proportional to
$|a_\e|^2$ for the integrated defect.}
Indeed,
\begin{equation}
 |x-a_\e|\geqslant|x|-|a_\e|>1>\delta_0
 \qquad(x\in K),
\end{equation}
hence $K\subset\Dext_{\e,\delta_0}$.
Restricting the preceding estimate \eqref{eq:tempmathcalB} to $K$ and using
the nonnegativity of $\Bp$, we obtain
\begin{equation}\label{eq:newmathcalB}
 \int_K\Bp(\bxi_0,\bxi_{a_\e})\dd
 \leqslant\Herr_\e(\delta_0)\longrightarrow0
 \qquad\text{as }\e\downarrow0.
\end{equation}
We next show that this vanishing defect forces
$|a_\e|\to0$.
We now bound this integral from below in terms of
$|a_\e|^2$. The formula
$\bxi_a(x)=(x-a)^\perp/|x-a|^2$ gives
\begin{equation}\label{eq:shift}
 |\bxi_0(x)-\bxi_{a_\e}(x)|
 =\frac{|a_\e|}{|x|\,|x-a_\e|},
 \qquad x\in\R^2\setminus\{0,a_\e\}.
\end{equation}
Since $|x|<3$ and $|x-a_\e|<4$ on $K$, it follows that
\begin{equation}
 |\bxi_0(x)-\bxi_{a_\e}(x)|
 \overset{\eqref{eq:shift}}{\geqslant}\frac{|a_\e|}{12},
 \qquad
 |\bxi_{a_\e}(x)|\geqslant\frac14,
 \qquad x\in K.
\end{equation}
Now we invoke the lower bound in \eqref{eq:convex}. Its quadratic term, integrated
over $K$, yields a constant $\kappa_p>0$,
depending only on $p$, such that
\begin{equation}\label{eq:center}
 \kappa_p|a_\e|^2
 \leqslant\int_K\Bp(\bxi_0,\bxi_{a_\e})\dd
 \leqslant\Herr_\e(\delta_0)\longrightarrow0.
\end{equation}
Thus $a_\e\to0$ for the original family, without
extracting a subsequence.

\smallskip\noindent
\emph{Step 2: the modulus tends uniformly to one away from the origin.}
Fix $r\in(0,1)$ and recall that
$A_r=B_1\setminus\overline{B_r}$.
We claim that
\begin{equation}\label{eq:mod_conv}
 \sup_{x\in\overline{A_r}}|\rh_\e(x)-1|
 \longrightarrow0
 \qquad\text{as }\e\downarrow0,
\end{equation}
where $\rh_\e=|\bU_\e|\in W^{1,p}_{\mathrm{loc}}(\R^2;\R)$.
We prove this by showing that every fixed loss of modulus
is confined to $B_r$ for sufficiently small $\e$.

Fix $\eta\in(0,1)$. Step~1 gives $a_\e\to0$, and
$C_{p,\eta}\lc\to0$ because $C_{p,\eta}$ is independent
of $\e$. Hence, for sufficiently small $\e$,
\begin{equation}
 \e<\e_{p,\eta},
 \qquad
 |a_\e|+C_{p,\eta}\lc<r.
\end{equation}
Proposition~\ref{prop:conc} then gives
\begin{equation}
 \{x\in B_1:|\bu_\e(x)|\leqslant1-\eta\}
 \subset B_{C_{p,\eta}\lc}(a_\e)
 \subset B_r.
\end{equation}
Since $\rh_\e=|\bu_\e|$ on $B_1$, it follows that
$\rh_\e(x)>1-\eta$ whenever $r\leqslant|x|<1$.
On $\partial B_1$, the prescribed boundary values give
$\rh_\e=1$. Together with $\rh_\e\leqslant1$, this yields
\begin{equation}
 \sup_{x\in\overline{A_r}}|\rh_\e(x)-1|
 \leqslant\eta
\end{equation}
for all sufficiently small $\e$. Since $\eta\in(0,1)$
was arbitrary, this proves \eqref{eq:mod_conv}.

\smallskip\noindent
\emph{Step 3: the phase defect and potential vanish on $A_r$.}
We claim that, for the fixed annulus
$A_r=B_1\setminus\overline{B_r}$,
\begin{equation}\label{eq:ann_conv}
 \|\nabla\rh_\e\|_{L^p(A_r)}\longrightarrow0,
 \qquad
 \|\bxi_\e-\bxi_0\|_{L^p(A_r)}\longrightarrow0,
 \qquad
 \frac1{4\e^2}\int_{A_r}(1-|\bu_\e|^2)^2\dd
 \longrightarrow0.
\end{equation}
Here $\rh_\e=|\bU_\e|$, and $\bxi_\e$ is the phase gradient
defined in \eqref{eq:phase}, which is well defined
on $A_r$ for small $\e$ by \eqref{eq:mod_conv}.
We first control the modulus gradient, the potential,
and a weighted phase difference using the defect estimate.
We then remove the weights and use $a_\e\to0$ to obtain
the phase convergence in \eqref{eq:ann_conv}.

Choose $\delta\in(0,\min\{r/2,r_*\})$, independently of $\e$.
By Step~1, $|a_\e|<r/2$ for sufficiently small $\e$.
Thus $|x-a_\e|>r/2>\delta$ for every
$x\in\overline{A_r}$, and hence
\begin{equation}
 \overline{A_r}\subset \Dext_{\e,\delta}
 =\R^2\setminus\overline{B_\delta(a_\e)}.
\end{equation}
Proposition~\ref{prop:conc}, with threshold $\eta=1/2$,
gives $|\bU_\e|\geqslant1/2$ on $\Dext_{\e,\delta}$ once
$\e<\e_{p,1/2}$ and $C_{p,1/2}\lc<\delta$.
Taking $\e$ smaller also ensures $\e<\e_*$ and
$M\lc\leqslant\delta$. Therefore both
Propositions~\ref{prop:excess} and \ref{prop:defect}
apply at this fixed radius.

Define the weighted phase difference
\begin{equation}\label{eq:phase_err}
 \bZ_\e=\rh_\e\bxi_\e
       -\rh_\e^{-1/(p-1)}\bxi_{a_\e}
 \in L^p(A_r;\R^2).
\end{equation}
Combining the two propositions, restricting the
nonnegative terms in \eqref{eq:defect} to $A_r$,
and applying \eqref{eq:convex}, we obtain
\begin{equation}\label{eq:grad_def}
 \frac1p\|\nabla\rh_\e\|_{L^p(A_r)}^p
 +b_p\|\bZ_\e\|_{L^p(A_r)}^p
 +\frac1{8\e^2}\int_{A_r}(1-\rh_\e^2)^2\dd
 \leqslant \Herr_\e(\delta)\longrightarrow0,
\end{equation}
where $\Herr_\e$ is defined in \eqref{eq:err}.
Each term tends to zero. Since $\rh_\e=|\bu_\e|$ on
$A_r$, this proves the first and third claims in
\eqref{eq:ann_conv}; the potential energy is
twice the last term in \eqref{eq:grad_def}.

To prove the remaining claim, we compare
$\bxi_{a_\e}$ with $\bxi_0$. The bound
$|x-a_\e|>r/2$ on $A_r$ and the identity in
\eqref{eq:shift} give
\begin{equation}
 \|\bxi_{a_\e}\|_{L^\infty(A_r)}\leqslant\frac2r,
 \qquad
 \|\bxi_{a_\e}-\bxi_0\|_{L^\infty(A_r)}
 \leqslant\frac{2|a_\e|}{r^2}\longrightarrow0.
\end{equation}
Now write
\begin{equation}\label{eq:unweight}
 \rh_\e(\bxi_\e-\bxi_0)
 =\bZ_\e
 +\bigl(\rh_\e^{-1/(p-1)}-\rh_\e\bigr)\bxi_{a_\e}
 +\rh_\e(\bxi_{a_\e}-\bxi_0).
\end{equation}
The first term tends to zero in $L^p(A_r;\R^2)$
by \eqref{eq:grad_def}. The second tends uniformly
to zero by \eqref{eq:mod_conv} and the uniform
bound on $\bxi_{a_\e}$. The third tends uniformly to zero
because $a_\e\to0$ and $\rh_\e\leqslant1$.
Thus $\rh_\e(\bxi_\e-\bxi_0)\to0$ in $L^p(A_r;\R^2)$.
Finally, $\rh_\e\geqslant1/2$ on $A_r$ for small $\e$,
so the same convergence holds without the factor
$\rh_\e$. This proves all three claims in
\eqref{eq:ann_conv}.

\smallskip\noindent
\emph{Step 4: convergence of the maps and their energies.}
We claim that
\begin{equation}
 \bu_\e\longrightarrow \bn
 \quad\text{strongly in }W^{1,p}(A_r;\C)
 \quad\text{and uniformly on }\overline{A_r},
 \qquad
 \E_\e(\bu_\e;A_r)\longrightarrow \cp(r^{2-p}-1).
\end{equation}
To recover the maps from the gradient estimates of Step~3,
we compare $\bu_\e$ with $\bn$ by complex multiplication.
The prescribed boundary values will fix the constant
phase that gradient convergence alone cannot determine.

Define
$\bz_\e=\bu_\e\overline{\bn}\in W^{1,p}(A_r;\C)$.
Since $(\bu_\e)_{|\partial B_1}=\bn$, its trace satisfies
$(\bz_\e)_{|\partial B_1}=1$.
Locally write $\bu_\e=\rh_\e e^{i\vp_\e}$ and
$\bn=e^{i\theta}$, so that
$\bz_\e=\rh_\e e^{i(\vp_\e-\theta)}$.
Using $\nabla\vp_\e=\bxi_\e$ and
$\nabla\theta=\bxi_0$, we obtain
\begin{equation}\label{eq:quotient}
 |\nabla \bz_\e|^2
 =|\nabla\rh_\e|^2
  +\rh_\e^2|\bxi_\e-\bxi_0|^2
 \qquad\text{almost everywhere on }A_r.
\end{equation}
Although the phases are chosen locally, every quantity
in this identity is globally defined. Since $\rh_\e\leqslant1$,
it follows that
\begin{equation}
 \|\nabla \bz_\e\|_{L^p(A_r)}^p
 \leqslant 2^{p/2-1}
 \left(
 \|\nabla\rh_\e\|_{L^p(A_r)}^p
 +\|\bxi_\e-\bxi_0\|_{L^p(A_r)}^p
 \right)\longrightarrow0
\end{equation}
by \eqref{eq:ann_conv}.

The function $\bz_\e-1$ has zero trace on $\partial B_1$.
Poincar\'e's inequality on $A_r$, with this outer boundary
condition, therefore gives a constant $C_{p,r}>0$ such that
\begin{equation}\label{eq:quot_conv}
 \|\bz_\e-1\|_{W^{1,p}(A_r)}
 \leqslant C_{p,r}\|\nabla \bz_\e\|_{L^p(A_r)}
 \longrightarrow0.
\end{equation}
Because $\bn\in W^{1,\infty}(A_r;\C)$ and
$\bu_\e-\bn=\bn(\bz_\e-1)$, multiplication by $\bn$ yields
$\bu_\e\to \bn$ strongly in $W^{1,p}(A_r;\C)$.
The continuous embedding
$W^{1,p}(A_r;\C)\hookrightarrow C^0(\overline{A_r};\C)$,
valid for $p>2$, then gives uniform convergence.
This proves \eqref{eq:conv}.

Finally, strong convergence of the gradients gives
\begin{equation}\label{eq:grad_lim}
 \frac1p\int_{A_r}|\nabla \bu_\e|^p\dd
 \longrightarrow
 \frac1p\int_{A_r}|\nabla \bn|^p\dd
 =\cp(r^{2-p}-1),
\end{equation}
where the last equality follows by polar integration.
Adding the potential energy, which tends to zero by
\eqref{eq:ann_conv}, proves \eqref{eq:en_lim}.
Since $r\in(0,1)$ was arbitrary and no subsequence was
extracted, the conclusions hold on every fixed annulus
for the original family of minimizers.
\end{proof}

\section*{Acknowledgements}
\textsc{G.D.F.} is a member of GNAMPA--INdAM. He gratefully
acknowledges partial financial support from the GNAMPA Project
CUP\_E53C25002010001, and from the University of Naples Federico
II through the FRA Project-B ``VarMoCry'' on 
\emph{Variational Analysis and Modeling of Liquid Crystals}.
Further support is acknowledged from the Italian Ministry of 
University and Research through the PRIN 2022 project 
\emph{Variational Analysis of Complex Systems in Material Science, Physics and Biology} (No.~2022HKBF5C).

 The work of AZ has been partially supported by the Basque 
 Government through the BERC 2022-2025 program and by 
 the Spanish State Research Agency through BCAM Severo 
 Ochoa excellence accreditation Severo Ochoa Grant CEX2021-001142-S funded by MICIU/AEI/10.13039/501100011033 and through project PID2023-146764NB-I00
funded by MICIU/AEI/10.13039/501100011033. 
AZ was also partially supported by a grant of 
the Ministry of Research, Innovation 
and Digitization, CNCS-UEFISCDI, project 
number PN-IV-P2-2.1-T-TE-2023-0219, within PNCDI IV.

\subsection*{Declarations} 
\emph{AI Usage:} Partially AI-assisted: the core mathematical approach was human-led, with Claude Opus 5.5 used to assist in drafting ideas or lemmas. The authors take full responsibility for the paper's contents. 
\emph{Data Availability:} The manuscript contains no associated data. 
\emph{Conflict of Interest:} The authors declare that they have no conflicts of interest. 
\emph{Author Contributions:} All authors contributed equally to this work.

\raggedbottom

\bibliographystyle{amsplain}
\bibliography{references}
\end{document}